\documentclass[journal]{IEEEtran}
\usepackage[usenames,dvipsnames]{color}
\usepackage{latexsym,psfig}
\usepackage{amsfonts}
\usepackage{amsmath,amssymb}
\usepackage{amsmath,epsfig}
\RequirePackage{color}
\usepackage{amssymb}
\usepackage{amsthm}
\usepackage{amsfonts}
\usepackage{mathrsfs}
\usepackage{algorithm}
\usepackage{algpseudocode}
\usepackage{times,lastpage,bm,xspace,booktabs}
\usepackage[center]{subfigure}
\usepackage{bibunits}
\usepackage{caption}
\usepackage{url}

\usepackage{mathrsfs}
\usepackage{amsmath,amssymb,latexsym}
\usepackage{amsfonts}
\usepackage{times,lastpage,bm,xspace,subfigure,bbm}
\usepackage{bibunits}
\usepackage{mathrsfs}
\usepackage{algorithm}

\def\Re{\mathbbm{R}}
\def\ones{\mathbbm{1}}
\newcommand{\ints}{{\mathbbm{Z}}}
\def \Poisson{{{\rm Poisson}}}
\def\deq{\triangleq}
\def\hf{\widehat{f}}

\def\ftrue{f^{\star}}

\def\argmin{\mathop{\arg \min}}
\def\argmax{\mathop{\arg \max}}
\def\minimize{\mathop{\textrm{minimize}}}
\def\st{\mathop{\textrm{\!subject to}}}
\def\grad{\nabla}

\def\pen{{\rm pen}}
\def\sP{\mathcal P}
\def\th{^{\rm th}}

\def\Int{\ensuremath{\mathbbm{Z}}}  
\def\ftrue{\ensuremath{f^\star}}
\def\thetatrue{\ensuremath{\theta^\star}}

\def\tif{\ensuremath{\widetilde{f}}}
\def\tv{\ensuremath{{\rm TV}}}
\newtheorem{theorem}{Theorem} 
\newtheorem{lemma}{Lemma}

\newenvironment{squishlist}
{   \begin{list}{$\bullet$}
    { \setlength{\itemsep}{2pt}      \setlength{\parsep}{2pt}
      \setlength{\topsep}{0pt}       \setlength{\partopsep}{0pt}
      \setlength{\leftmargin}{1.5em} \setlength{\labelwidth}{1em}
      \setlength{\labelsep}{0.5em} } }
      {\end{list}}

\usepackage{booktabs}

\newcommand{\spiral}{SPIRAL}
\usepackage[normalem]{ulem} 

\begin{document}

\title{This is \spiral-TAP: Sparse Poisson Intensity Reconstruction
  ALgorithms -- Theory and Practice}
\author{Zachary~T.~Harmany,
	~\IEEEmembership{Student Member,~IEEE,}
        Roummel~F.~Marcia,
        ~\IEEEmembership{Member,~IEEE,}\\
        Rebecca~M.~Willett
        ,~\IEEEmembership{Member,~IEEE}        %
\thanks{Z.~T.~Harmany and R.~M.~Willett are with 
the Department of Electrical and Computer Engineering,
Duke University, Durham, NC 27708, \textbf{USA}.
R.~F.~Marcia is with 
the School of Natural Sciences,
The University of California, Merced, Merced, CA 95343, \textbf{USA}}}
%
%
\markboth{Submitted to IEEE Transactions on Image Processing}%
{Shell \MakeLowercase{\textit{et al.}}: Bare Demo of IEEEtran.cls for Journals}
%
\maketitle

\begin{abstract}
  \boldmath The observations in many applications consist of counts of
  discrete events, such as photons hitting a detector, which cannot be
  effectively modeled using an additive bounded or Gaussian noise
  model, and instead require a Poisson noise model. As a result,
  accurate reconstruction of a spatially or temporally distributed
  phenomenon ($\ftrue$) from Poisson data ($y$) cannot be effectively
  accomplished by minimizing a conventional penalized least-squares
  objective function. The problem addressed in this paper is the
  estimation of $\ftrue$ from $y$ in an inverse problem setting, where
  (a) the number of unknowns may potentially be larger than the number
  of observations and (b) $\ftrue$ admits a sparse approximation.  The
  optimization formulation considered in this paper uses a penalized
  negative Poisson log-likelihood objective function with
  nonnegativity constraints (since Poisson intensities are naturally
  nonnegative).  In particular, the proposed approach incorporates key
  ideas of using separable quadratic approximations to the objective
  function at each iteration and penalization terms related to
  $\ell_1$ norms of coefficient vectors, total variation seminorms, and
  partition-based multiscale estimation methods.
\end{abstract}

%
\IEEEpeerreviewmaketitle

\section{Introduction} 

\IEEEPARstart{I}{n} a variety of applications, ranging from nuclear
medicine to night vision and from astronomy to traffic analysis, data
are collected by counting a series of discrete events, such as photons
hitting a detector or vehicles passing a sensor.  The measurements are
often inherently noisy due to low count levels, and we wish to
reconstruct salient features of the underlying phenomenon from these
noisy measurements as accurately as possible.  The inhomogeneous
Poisson process model \cite{snyder} has been used effectively in many
such contexts. Under the Poisson assumption, we can write our
observation model as
\begin{equation}
y \sim \Poisson(A\ftrue),
\label{eq:PoissonModel}
\end{equation}
where $\ftrue \in \Re_+^n$ is the signal or image of interest, $A \in
\Re_+^{m \times n}$ linearly projects the scene onto an
$m$-dimensional set of observations, and $y \in \Int_+^m$ is a
length-$m$ vector of observed photon counts.

The problem addressed in this paper is the estimation of $\ftrue$ from
$y$ when $\ftrue$ is sparse or compressible in some basis $W$ (i.e.,
$\ftrue = W \thetatrue$ and $\thetatrue$ admits an accurate sparse
approximation) and the number of unknowns $n$ may be larger than the
number of observations $m$. This challenging problem has clear
connections to \emph{compressed sensing}
(CS) \cite{RIP,CS:donoho,CS:candes1,CS:candes2}, but arises in a number
of other settings as well, such as tomographic reconstruction in
nuclear medicine, superresolution image reconstruction in astronomy,
and deblurring in confocal microscopy.  In recent work \cite{WilRag09,PoissonCS_tsp},
 we explored some of the theoretical challenges
associated with CS in a Poisson noise setting, and in particular
highlighted two key differences between the conventional CS problem
and the Poisson CS problem:
\begin{squishlist}
\item unlike many sensing matrices in the CS literature, the matrix
  $A$ must contain all nonnegative elements, and
\item the intensity $\ftrue$ is nonnegative, and hence any estimate $\hf$ thereof must also be nonnegative.
\end{squishlist}
The nonnegativity of $\hf$
and $A$ results in challenging optimization
problems. In particular, the restriction that $\hf$ is nonnegative
introduces a set of inequality constraints into the minimization
setup; as shown in \cite{fesslerPoissonCS, harmanyPCS}, these
constraints are simple to satisfy when $\hf$ is sparse in the canonical
basis, but they introduce significant challenges when enforcing
sparsity in an arbitrary basis. 

\subsection{Problem formulation}
Under the Poisson model \eqref{eq:PoissonModel}, the probability of
observing a particular vector of counts $y$ is given by
\begin{equation}
p \left(y \middle| A \ftrue \right)
= \prod_{i=1}^m \frac{(e_i^TA\ftrue)^{y_i}}{y_i!} \exp\left(-e_i^TA\ftrue\right)
\label{eq:PoissonLikelihood}
\end{equation}
where $e_i$ is the $i\th$ canonical basis unit vector.
Thus the negative Poisson log-likelihood is given by
\begin{equation}
F(f) = \ones^TAf - \sum_{i=1}^m y_i \log(e_i^T A f),
\label{eq:negpoisson}
\end{equation}
where $\ones$ is an $m$-vector of ones.
(Here we neglect $\log(y_i!)$ terms since they are constant with
respect to $f$ and hence do not impact our objective.) We will see
later that in order to avoid the singularity at $f=0$, it is
advantageous to introduce a small parameter $\beta > 0$ where 
typically $\beta \ll 1$: 
\begin{equation}
F(f) \equiv \ones^TAf - \sum_{i=1}^m y_i \log(e_i^T A f + \beta).
\label{eq:negpoisson_epsilon}
\end{equation}
Similar techniques are used in \cite{Fessler98aparaboloidal} in the form of
known background intensities, where the observations are modeled by $y \sim \Poisson (A\ftrue + b)$ with $b \in \Re_+^n$ known a priori. If $\min(b)$ is large enough $\beta$ is not explicitly needed, however in our experiments we assume $b = 0$ and use $\beta = 1\times 10^{-10}$, much smaller than any reasonable background count magnitude.
This small parameter will also appear in the gradient
\begin{equation}
\nabla F(f) = A^T\ones-\sum_{i=1}^m \frac{y_i}{e_i^T A f + \beta} A^T e_i,
\label{eq:gradient}
\end{equation}
and the Hessian
\begin{equation}
\nabla^2 F(f) = A^T \left[\sum_{i=1}^m \frac{y_i}{(e_i^T A f + \beta)^2}e_i e_i^T\right] A.
\label{eq:hessian}
\end{equation}

Our Poisson reconstruction algorithms take the form of the following constrained optimization problem:
\begin{equation}
\begin{aligned}
&\underset{f \in \mathbbm{R}^n}{\textrm{minimize}} \ \Phi(f) \equiv F(f) + \tau \pen(f) \\
&\textrm{subject to} \ f \geq 0 ,
\label{eq:spiral_goal}
\end{aligned}
\end{equation}
where
\begin{itemize}
\item $F:\Re^n \to \Re$ is the negative Poisson log-likelihood in
	\eqref{eq:negpoisson_epsilon}, and
\item $\pen : \Re^n \to \Re$ is a finite, usually nonsmooth, and
  potentially nonconvex penalty functional.
\end{itemize}
In this paper, we will consider several variants of the penalty term, including
$\|f\|_1$, $\|W^T\! \! f\|_1$ for some arbitrary orthonormal basis
$W$, a total-variation seminorm $\|f\|_\tv$, and a complexity
regularization penalty based upon recursive dyadic partitions. We
refer to our approach as \spiral\ (Sparse Poisson Intensity
Reconstruction ALgorithm). 

\subsection{Related work}
Various regularization techniques are often employed to compensate for
the ill-posedness of the estimation problem. Outside the Poisson
context, for example in the presence of additive white Gaussian noise,
regularization methods based on wavelet or curvelet sparsity \cite{CS:candes2},
models of wavelets' clustering and persistence
properties \cite{baraniukModelBasedCS}, and a variety of other
penalties (cf. \cite{sparsa}) have proven successful.

The key challenge in Poisson intensity estimation problems is that the
mean and the variance of the observed counts are the same. As a
result, conventional approaches based on a penalized least-squares
objective function (cf. \cite{sparsa,gpsr}) will yield suboptimal
results when applied to Poisson data with low intensities.
Variance-stabilizing transforms (VSTs), such as the Anscombe transform
\cite{anscombe} and the Haar-Fisz \cite{fryzlewicz2008data,HaarFisz}
transform, are widely used to address this issue and to approximate
the Poisson observations by Gaussian random variables \cite{
  donoho93nonlinear,Poiss_Anscombe}. Jansen proposes a wavelet based
Poisson estimation method based on data-adaptive VSTs and Bayesian
priors on the stabilized coefficients
\cite{jansen2006multiscale}. However, as pointed out in
\cite{kola99sinica,willett:Review}, such approximations are inaccurate
when the observed number of photons per pixel or voxel is very low and
tend to oversmooth the resulting estimate. In a more recent work,
Zhang et al. \cite{Starck_Poisson} propose a multiscale
variance-stabilizing transform (MSVST) which applies a VST to the
empirical wavelet, ridgelet or curvelet transform coefficients.
However, theoretical analysis of this approach is not
available and it is not clear how to extend the MSVST to Poisson
inverse problems.

Several authors have investigated reconstruction algorithms
specifically designed for Poisson noise without the need for VSTs. In
\cite{nowkol:00,multiscaleMAP}, Nowak and Kolaczyk describe multiscale
Bayesian approach in an Expectation-Maximization (EM) framework to
reconstruct Poisson intensities.  In their seminal paper \cite{kola2},
Kolaczyk and Nowak present a multiscale framework for likelihoods
similar to wavelet analysis and propose a denoising algorithm based on
the complexity-penalized likelihood estimation (CPLE). The CPLE
objective function is minimax optimal over a wide range of isotropic
likelihood models. There are several variants of the CPLE method
depending upon the nature of the image or signal being denoised
\cite{willett:density,glm,willett:tmi03}.

Regularization based on a {\em total variation (TV)} seminorm has also
garnered significant recent attention (cf.
\cite{chanTV,BeckTV}). This seminorm is described in detail below;
in general, it measures how much an image varies across pixels, so
that a highly textured or noisy image will have a large TV seminorm,
while a smooth or piecewise constant image would have a relatively
small TV seminorm. This is often a useful alternative to wavelet-based
regularizers, which are also designed to be small for piecewise smooth
images but can result in spurious large, isolated wavelet coefficients
and related image artifacts.

In the context of Poisson inverse problems, however, adaptation of
these regularization methods can be challenging for two main 
reasons. First, the negative Poisson log-likelihood used in the
formulation of an objective function often requires the application of
relatively sophisticated optimization theory principles. Second,
because Poisson intensities are inherently nonnegative, the resulting
optimization problem must be solved over a constrained feasible set,
increasing the complexity of most algorithms.  Some recent headway has
been made using multiscale or smoothness-based penalties
\cite{willett:density,fesslerSPSOS,spieral}.
In one recent work \cite{starckPoissonCS}, the Poisson statistical
model is bypassed in favor of an additive Gaussian noise model through
the use of the Anscombe variance stabilizing transform. This
statistical simplification is not without cost, as the linear
projections of the scene must now be characterized as nonlinear
observations.  Other recent efforts \cite{marioPCS,setzer} solve
Poisson image reconstruction problems with TV seminorm regularization
using a split Bregman approach \cite{GoldsteinOsher}, but the proposed
methods involves a matrix inverse operation which can be extremely
difficult to compute for large problems outside of deconvolution
settings. TV seminorm regularization is also explored in the context of a
Richardson-Lucy algorithm \cite{zerubiaPoissonTV}, but nonnegativity
and convergence were not explicitly addressed.  Finally, the
approaches in \cite{combettes,chaux} apply proximal functions to solve
more general constrained convex minimization problems.  These methods
use projection to obtain feasible iterates (i.e., nonnegative
intensity values), which may be difficult for recovering signals that
are sparse in a noncanonical basis.

\subsection{Contributions of the proposed method}

In the proposed work, we present a general algorithmic framework for solving Poisson inverse problems. This framework requires no special structure in the sensing matrix $A$, yet can take advantage fast matrix-vector multiplications when available. The success of this approach hinges on the ability to solve certain constrained subproblems quickly. We show that this can be done for a variety of regularization schemes. Additionally, we consider a step selection procedure that is better suited for the Poisson log-likelihood versus the corresponding generalizations from the Gaussian log-likelihood. In the subsequent analysis, we establish global convergence of the constrained optimization under a set of mild assumptions which are easily satisfied in practice. We establish properties of the solution set, which yield conditions for a unique solution to the minimization. Our approach is then supported by numerical simulations which show state-of-the-art performance on a simulated limited-angle emission tomography inverse problem.

\section{Algorithms}
\label{sec:algorithms}

Our approach to solve the minimization problem \eqref{eq:spiral_goal}
employs sequential quadratic approximations to the Poisson
log-likelihood $F(f)$. More specifically, at iteration $k$ we compute
a separable quadratic approximation to $F(f)$ using its second-order
Taylor series approximation at $f^k$; this approximation is denoted
$F^k(f)$.  The next iterate is then given by
\begin{equation}\begin{aligned}
f^{k+1} =\ &\underset{f \in \mathbbm{R}^n}{\arg \min} \ F^k(f) + \tau \pen(f) \\
	& \textrm{subject to} \ f \geq 0.
\label{eq:spiral_sub_general}
\end{aligned}\end{equation}
Similar to the framework described in \cite{sparsa}, $F^k$ is
a second-order Taylor series approximation to $F$ with the Hessian
$\nabla^2 F(f^k)$ approximated by a scaled identity matrix $\alpha_k
I$, with $\alpha_k > 0$. This yields
\begin{equation}
F^k(f) = F(f^k) + (f-f^k)^T\nabla F(f^k) + \tfrac{\alpha_k}{2}\|f - f^k\|_2^2.
\label{eq:separable_approx}
\end{equation}
With this separable approximation, simple manipulation of 
\eqref{eq:spiral_sub_general} yields a sequence of subproblems of the form
\begin{equation}
\begin{aligned}
f^{k+1} =\ &\underset{f \in \mathbbm{R}^n}{\arg \min} \ \phi^k(f) = \tfrac{1}{2}\|f - s^k\|_2^2 + \tfrac{\tau}{\alpha_k} \pen(f) \\
&\textrm{subject to} \ f \geq 0,
\label{eq:spiral_sub}
\end{aligned}
\end{equation}
where 
\begin{equation}	\label{eq:sk}
s^k = f^k - \tfrac{1}{\alpha_k} \nabla F(f^{k}).
\end{equation}
This formulation has the benefit of being easily recognized as a nonnegatively constrained
$\ell_2$ denoising of the gradient descent step $s^k$. 

The parameter $\alpha_k$ is chosen via a sequence of two repeated steps.  First, a modified Barzilai-Borwein (BB) method \cite{barzilai} is used to choose the initial value of $\alpha_k$.  With $\delta^k = f^k - f^{k-1}$, we initially choose
\begin{equation}
\alpha_k = \frac{(\delta^k)^T \nabla^2 F(f^k) \delta^k}{\|\delta^k\|_2^2}
= \frac{\| \sqrt{y} \cdot (A\delta^k) / (Af^k + \beta) \|_2^2}{\| \delta^k\|_2^2},
\label{eq:alphainit}
\end{equation}
with $\alpha_k$ safeguarded to be within the range $[\alpha_\text{min}, \alpha_\text{max}]$.  
Here $\sqrt{\cdot}$, $\cdot$, and $/$ are understood as component-wise. This method allows our
separable approximation \eqref{eq:separable_approx} to capture the curvature of the Poisson log-likelihood $F(f)$ along the most recent
step $\delta^k$, in the vicinity of the current iterate $f^k$.

The astute reader will recognize \eqref{eq:alphainit} as a Rayleigh quotient, and hence always has a value within the spectrum of the Hessian. Since $F$ is convex, this guarantees $\alpha_k \ge 0$, even without the safeguard. The traditional BB scheme chooses $\alpha_k = (\gamma^k)^T \delta^k/\|\delta^k\|_2^2$, where $\gamma^k = \grad F (f^k) - \grad F (f^{k-1})$. When applied to a quadratic objective, such as arising from considering the Gaussian log-likelihood, these two schemes are identical. However, when applied to the Poisson log-likelihood, our choice is more strongly effected by the curvature at the current iteration.  Our modified BB choice is no more expensive to compute than the original BB scheme since $A\delta^k$ and $Af^k$ are already available from the gradient computations, and computing $\alpha_k$ in \eqref{eq:alphainit} is a simple sequence of $O(n)$ operations.

This initial choice of $\alpha_k$ is used if the resulting solution of \eqref{eq:spiral_sub} satisfies the acceptance criteria
\begin{equation}
\Phi(f^{k+1}) \leq \max_{i = [k-M]_+,\dots,k} \Phi(f^i) - \tfrac{\sigma \alpha_k}{2} \|f^{k+1} - f^{k} \|_2^2,
\label{eq:spiral_accept}
\end{equation}
where $M$ is a nonnegative integer, $\sigma \in (0,1)$ is a small
constant, and the operation 
$[ \, \cdot \, ]_+ = \max \{0, \, \cdot \, \}$.
If it fails this acceptance criteria, $\alpha_k$ is
repeatedly increased by a factor $\eta$ until the solution to
\eqref{eq:spiral_sub} satisfies \eqref{eq:spiral_accept}.  This gentle
criteria allows the nonmonotonic objective behavior characteristic of
the Barzilai-Borwein methods \cite{sparsa,barzilai}, yet enforces that
the next iterate have a slightly smaller objective than the largest
value over the past $M$
iterations.  Note that choosing $M=0$ results in a purely monotonic algorithm.
For clarity, we describe our general procedure in Algorithm
\ref{spiral_alg}.

\begin{algorithm} 
\caption{Sparse Poisson Intensity \\ Reconstruction ALgorithm (\spiral)}
\label{spiral_alg} 
\begin{algorithmic}[1] 
\State \textbf{Initialize}  Choose $\eta > 1$, $\sigma \in (0,1)$, $M \in \mathbbm{Z}_+$, $0 < \alpha_\text{min} \leq \alpha_\text{max}$, and initial solution $f^0$.  Start iteration counter $k \gets 0$.
\Repeat
\State choose $\alpha_k \in [\alpha_\text{min}, \alpha_\text{max}]$ by \eqref{eq:alphainit}
\State $f^{k+1} \gets$ solution of \eqref{eq:spiral_sub}
\While{$f^{k+1}$ does not satisfy \eqref{eq:spiral_accept}}
\State $\alpha_k \gets \eta \alpha_k$
\State $f^{k+1} \gets$ solution of \eqref{eq:spiral_sub}
\EndWhile
\State $k \gets k + 1$
\Until{stopping criterion is satisfied}
\end{algorithmic} 
\end{algorithm}

\subsection{Canonical basis with sparsity penalty}

\label{sec:canonical}
When  
$\pen(f) = \| f \|_1$,
the minimization subproblem
(\ref{eq:spiral_sub}) 
has the following analytic solution:
\begin{equation*}
	f^{k+1} = \left [ s^k - \tfrac{\tau}{\alpha_k} \ones \right ]_+,
\end{equation*}
with $[ \, \cdot \, ]_+$ acting component-wise.
Thus solving (\ref{eq:spiral_goal}) subject to nonnegativity constraints with an $\ell_1$ penalty function measuring sparsity in the canonical basis is straightforward.
An alternative algorithm for solving this Poisson inverse
problem with sparsity in the canonical basis was also explored in the
recent literature \cite{fesslerPoissonCS}.

\subsection{Non-canonical basis} 

Now suppose that the signal of interest is sparse in some other basis.  Then
the $\ell_1$ penalty term is given by
\begin{equation*}
\pen(f) \deq \|W^T\!\!f\|_1 =  \|\theta\|_1,
\end{equation*}
where 
\begin{equation}
	\theta \deq W^T\!f	\label{eq:Wtheta}
\end{equation} 
for some orthonormal basis $W$.
 When the reconstruction $\hat{f} = W\hat{\theta}$ must be nonnegative (i.e.,
$W\hat{\theta} \ge 0$), the minimization problem
\begin{align}	\nonumber
	\theta^{k+1} \deq \ &
	\underset{\theta \in \Re^n}{\arg \min } \quad
	\phi^{k}(\theta) \deq \tfrac{1}{2} \| \theta - s^k \|_2^2 + \tfrac{\tau}{\alpha_k} \| \theta \|_1, \\
	& \textrm{subject to \ } W\theta \ge 0
	\label{eq:sparsa_con} 
\end{align}
no longer has an analytic solution necessarily.  
We can solve this minimization problem by solving its Lagrangian dual.
First, we reformulate 
(\ref{eq:sparsa_con}) so that its objective function $\phi^k(\theta)$ is differentiable 
by defining $u, v \ge 0$ such that $\theta = u - v$.  The minimization problem (\ref{eq:sparsa_con})
becomes
\begin{align}	
	\nonumber
	(u^{k+1}, v^{k+1}) \deq \ &
	\underset{u,v \in \Re^{n}}{\arg \min } \quad
	\tfrac{1}{2} \| u -  v - s^k \|_2^2 + \tfrac{\tau}{\alpha_k} \ones^T\!(u+v) \\
	\label{eq:primal}
	& \textrm{subject to \ } u,v \ge 0, \quad W(u-v) \ge 0, \\[-.5cm] \nonumber
\end{align}
which has twice as many parameters and has additional nonnegativity constraints on the new parameters, 
but now has a differentiable objective function.
The Lagrangian function corresponding to (\ref{eq:primal})
is given by
\begin{align*}
	\mathscr{L}(u,v,\lambda_1, \lambda_2, \lambda_3)
	&=  \tfrac{1}{2} \| u-v - s^k \|_2^2 + \tfrac{\tau}{\alpha_k}\ones^T(u+v) \\
	& \qquad 	\ - \lambda_1^Tu - \lambda_2^Tv - \lambda_3^TW(u-v),
\end{align*}
where $\lambda_1, \lambda_2, \lambda_3 \in \Re^n$
are the Lagrange multipliers corresponding to the constraints in (\ref{eq:primal}).  
Setting the derivative of $\mathscr{L}$ with respect to $u$ and $v$ to zero, 
we obtain
\begin{align}
	\label{eq:uv} 
	u-v &= s^k  + \lambda_1 - \tfrac{\tau}{\alpha_k} \ones  + W^T\! \lambda_3, \text{ and} \\
	\lambda_2 &= \tfrac{2\tau}{\alpha_k} \ones - \lambda_1, \nonumber
\end{align}
which leads to the Lagrangian dual function
\begin{equation*}
	g(\lambda_1,\lambda_3) =  - \frac12\| s^k + \lambda_1 - \tfrac{\tau}{\alpha_k} \ones 
	+ W^T\! \lambda_3 \|_2^2 + \frac12 \| s^k \|_2^2.
\end{equation*}
We define $\gamma \deq \lambda_1 - \tfrac{\tau}{\alpha_k} \ones$.  For the Lagrange dual problem corresponding to
(\ref{eq:primal}), the Lagrange multipliers $\lambda_1, \lambda_2, \lambda_3 \ge 0$.
Since $\lambda_2 = \tfrac{2\tau}{\alpha_k} \ones - \lambda_1$,
then $- \tfrac{\tau}{\alpha_k} \ones \le \gamma \le \tfrac{\tau}{\alpha_k} \ones$.  Also, let $\lambda = \lambda_3$.
The Lagrange dual problem associated with this problem is thus given by
\begin{eqnarray}	\nonumber
	&& \underset{\lambda, \gamma \in \Re^n}{\textrm{minimize}} \ \ 
	h(\lambda, \gamma) \deq
	\tfrac{1}{2} \| s^k + \gamma  + W^T\! \lambda  \|_2^2 - \tfrac{1}{2} \| s^k \|_2^2  \quad \\ 
	&& \textrm{subject to \ } \lambda \ge 0, \ \ \  - \tfrac{\tau}{\alpha_k} \ones \le \gamma \le \tfrac{\tau}{\alpha_k} \ones 
	\label{eq:dual} \\[-.3cm] \nonumber
\end{eqnarray}
and at the optimal values $\gamma^{\star}$ and
$\lambda^{\star}$, the primal iterate $\theta^{k+1}$  is given by
$$
	\theta^{k+1} \deq u^{k+1} - v^{k+1} = s^k + \gamma^{\star} + W^T\! \lambda^{\star}.
$$
We note that the minimizers of the primal problem (\ref{eq:primal})
and its dual (\ref{eq:dual}) satisfy $\phi^k(\theta^{k+1}) = -h(\gamma^{\star}, \lambda^{\star})$
since (\ref{eq:primal}) satisfies (a weakened) Slater's condition \cite{convexoptimization}.
In addition, the function $-h(\gamma, \lambda)$ is a lower bound on $\phi^{k}(\theta)$
at any dual feasible point.

The objective function of (\ref{eq:dual}) can be minimized by alternatingly solving for $\lambda $ and $\gamma$,
which is accomplished by taking the partial derivatives of $h(\lambda, \gamma)$ and setting them to zero. 
Each component is then constrained to satisfy the bounds in (\ref{eq:dual}).  At the $j^{\rm th}$ iteration, 
the variables can, thus, be defined as follows:
\begin{eqnarray}	\nonumber
	\gamma^{(j)} \!\!\! &=& \!\! \text{mid}  \Big \{  -\tfrac{\tau}{\alpha_k} \ones, -s^k -W^T \! \lambda^{(j-1)},
			\tfrac{\tau}{\alpha_k} \ones \Big \} \\
	\lambda^{(j)} \!\!\! &=& \!\! \left [ -W \left (s^k + \gamma^{(j)} \right ) \right]_+,	\label{eq:lambdaj}
\end{eqnarray}
where the operator $\text{mid}\{a,b,c\}$ chooses the middle value of the three arguments component-wise.  
Note that at the end of each iteration $j$, the 
approximate solution $\theta^{(j)} \deq s^k + \gamma^{(j)} + W^T \lambda^{(j)}$ to (\ref{eq:sparsa_con}) 
is feasible with respect to the constraint $W\theta \ge 0$:
\begin{eqnarray}	\nonumber
	W\theta^{(j)} &=& 
	Ws^k + W\gamma^{(j)}  + \lambda^{(j)}  \\	\nonumber
	&=& W(s^k + \gamma^{(j)})  + \left [ -W(s^k + \gamma^{(j)})  \right ]_+  \\	
	&=& \left [ W \left (s^k + \gamma^{(j)} \right )  \right ]_+  \label{eq:Wthetaj}  \\
	&\ge& 0.	 \nonumber
\end{eqnarray}
Thus, we can terminate the iterations for the dual problem early
and still obtain a feasible point.  By solving the subproblem ``loosely"
(see Sec. \ref{sec:termination}), the algorithm can converge faster but to 
a potentially less accurate estimate.  In Sec.~\ref{sec:experiments}, we show
results from both ``loose" and ``tight" criteria for solving the subproblem.
We call this approach \spiral-$\ell_1$.

\subsection{Total variation penalty}
We also consider an approach where the TV seminorm is used as a penalization scheme, i.e., $\pen(f) = \|f\|_\tv$. We define the anisotropic TV seminorm as
\begin{equation*}\begin{aligned}
\|f\|_\tv  \deq \sum_{k=1}^{\sqrt{n}-1} \sum_{l=1}^{\sqrt{n}}
	|f_{k,l} - f_{k+1,l}|	
 + \sum_{k=1}^{\sqrt{n}} \sum_{l=1}^{\sqrt{n}-1} |f_{k,l} - f_{k,l+1}|,
\end{aligned}\end{equation*}
where -- for simplicity of presentation, not algorithmic necessity -- we use a slight abuse of notation by using 2D pixel indices instead of vector indices by assuming that $f \in \Re^n$ is a square $\sqrt{n} \times \sqrt{n}$ image. This highlights the fact that the TV seminorm is simply a measure of the magnitude of all vertical and horizontal first-order differences. Said differently, one can think of an image with a small TV seminorm as one that is sparse with respect to an overcomplete representation of these first-order differences (neglecting the mean of $f$), meaning the image has few abrupt changes in pixel intensity yet many regions of homogeneous signal level.  This property makes TV especially well suited for image denoising and inverse problems.

Recent work by Beck and Teboulle \cite{BeckTV} presents a fast computational method for
solving the TV-regularized problem
\begin{equation}
\begin{array}{rll}
	\tif \ = & \displaystyle \argmin_{f \in \Re^n}  
	& \displaystyle \frac{1}{2}\|\widetilde{A}f -b \|_2^2 + \lambda \| f \|_\tv \\
	& \displaystyle \st & f \in C,
\end{array} \label{eq:tvbeck} 
\end{equation}
where $\lambda>0$ is a tuning parameter, $C$ is a closed convex set and
$\widetilde{A}$ is a linear, {\em spatially invariant} blur operator. This method utilizes a gradient-based optimization approach founded on a monotone iterative shrinkage and thresholding algorithm.
When we choose $\tif = f^{k+1}$, $\widetilde{A} = I$, $b = s^k$, $\lambda = \tau/\alpha_k$, and $C = \{f \in \Re^n : f \geq 0\}$, \eqref{eq:tvbeck} then reduces to an $\ell_2$ denoising of $s^k$ with a total variation regularizer:
\begin{equation}
\begin{array}{rll}
	f^{k+1} \ = & \displaystyle \argmin_{f \in \Re^n}  
	& \displaystyle \tfrac{1}{2}\|f -s^k\|_2^2 + \tfrac{\tau}{\alpha_k} \| f \|_\tv \\
	& \displaystyle \st & f \geq 0,
\end{array} \label{eq:tvspiral} 
\end{equation}
precisely the form \eqref{eq:spiral_sub} required for use in our
algorithmic framework. We call this approach with a total variation penalty \spiral-TV.

\subsection{Partition-based methods}

\begin{figure}[htb]
	\centering
	\subfigure[]{\includegraphics[width=1.1in]{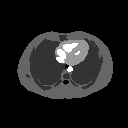}}
	\subfigure[]{\includegraphics[width=1.1in]{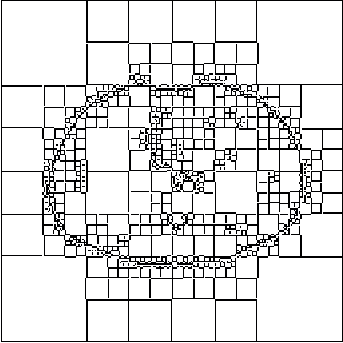}}
	\subfigure[]{\includegraphics[width=1.1in]{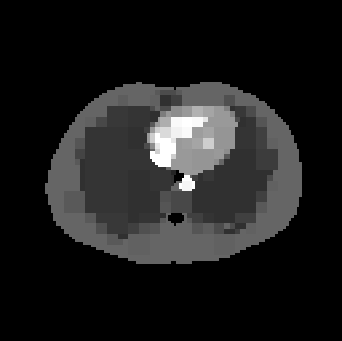}}
	\caption{	Example of a partition-based approximation. (a) True
          image. (b) Recursive dyadic partition (RDP). (c)
          RDP-based approximation of the true image.}
	\label{fig:Partition}
\end{figure}

An alternative to the $\ell_1$- and TV-norm penalties 
can be formulated using model-based estimates that utilize structure in the coefficients 
beyond that of a sparse representation.  In particular, we build upon the
framework of {\em recursive dyadic partitions
  (RDP)}, which we summarize here and are described in detail in
\cite{willett:density}. Let $\sP$ be the class of all recursive dyadic
partitions of $[0,1]^2$ where each cell in the partition has a
sidelength at least $1/\sqrt{n}$, and let $P \in \sP$ be a candidate
partition. The intensity on $P$, denoted $f(P)$, is calculated using a
nonnegative least-squares method to fit a model (such as a constant)
to $s^k$ in (\ref{eq:sparsa_con}) in each cell in the RDP.  As an
example, consider Fig.~\ref{fig:Partition}.  Here we approximate the
true image (Fig.~\ref{fig:Partition}(a)) on the recursive dyadic
partition defined in Fig.~\ref{fig:Partition}(b)).  The result is a
piecewise constant approximation to the emission image
(Fig.~\ref{fig:Partition}(c)).  We see that the partition model is
able to accurately capture the image in clear multiresolution fashion:
large homogeneous regions are well-modeled by large cells, whereas
edges are approximated via the deeper recursive partitioning.
Furthermore, a penalty can be assigned to the resulting
estimator which is proportional to $|P|$, the number of cells in $P$.
Thus we set
\begin{equation} \label{eq:partition}
\begin{array}{rcl}
\widehat{P}^{k+1} &=& \displaystyle \argmin_{P \in \sP} \ \tfrac{1}{2} \|s^k - f(P)\|_2^2 + \tfrac{\tau}{\alpha_k} |P|,
\medskip \\
f^{k+1} &=& f(\widehat{P}^{k+1}).
\end{array}
\end{equation}
A search over $\sP$ can be computed quickly using a dynamic
program. When using constant partition cell models, the nonnegative
least-squares fits can be computed non-iteratively in each cell by
simply using the maximum of the average of $s^k$ in that cell and
zero. Because of this, enforcing the constraints is trivial and can be
accomplished very quickly.  The disadvantage of using constant model
fits is that it yields piecewise constant estimates. However, a
cycle-spun translation-invariant (TI) version of this
approach \cite{willett:density} can be implemented with high
computational efficiency and can be used for solving this nonnegative
regularized least-squares subproblem that results in a much smoother
estimator.  We refer to these approaches as \spiral-RDP and
\spiral-RDP-TI.

It can be shown that partition-based denoising methods such as this
are closely related to Haar wavelet denoising with an important
hereditary constraint placed on the thresholded coefficients -- if a
parent coefficient is thresholded, then its children coefficients must
also be thresholded \cite{willett:density}. This constraint is akin to
wavelet-tree ideas which exploit persistence of significant wavelet
coefficients across scales and have recently been shown highly useful
in compressed sensing settings \cite{baraniukModelBasedCS}. Since in this context the penalty can be thought of an $\ell_0$ measure, the resulting RDP-based penalty function is not a convex function.  Hence we can only guarantee convergence to a local minimizer.  Given a sufficiently accurate initialization, the RDP-based method performs competitively in both speed and accuracy to convex penalties, for which global convergence is assured.

\section{Algorithmic details}

\subsection{Computational complexity and nonmonotonicity}

It should be noted that in determining an acceptable next
iterate $f^{k+1}$ in \eqref{eq:spiral_sub_general}, the main
computational burden is often checking the acceptance criteria
\eqref{eq:spiral_accept} as testing each candidate solution involves
recomputing $A f^{k+1}$. In the worst case, $A$ may be a dense
unstructured matrix for which computing matrix-vector products is
costly.  Although enforcing near-monotonicity of the objective aids
convergence to a more accurate solution to the reconstruction problem
\eqref{eq:spiral_goal}, significant demonstrable gains in
computational speed may be achieved by forgoing this criteria and
simply accepting the choice \eqref{eq:alphainit}.  In this case, an
efficient implementation of our algorithm only requires two
matrix-vector multiplications involving $A$ per iteration, one
computing $Af^k$ for defining $\alpha_k$ in \eqref{eq:alphainit}, and the other using $A^T$ in the computation of
$\nabla F(f^k)$ for defining the gradient descent result $s^k$ in \eqref{eq:sk}.  
Since  \eqref{eq:spiral_sub} requires only $\alpha_k$
and $s^k$, there are no matrix-vector multiplications involving $A$ 
in the denoising subproblem.

\subsection{Initialization}

While global convergence proofs -- such as the one in Section~\ref{sec:convergence} -- guarantee convergence for any initial point $f^0 \ge 0$, the choice of initialization is an important practical consideration.  Iterative algorithms are rarely allowed to execute long enough to converge to an optimal solution, meaning the approximate solution may be strongly dependent on the starting point of the algorithm.  Further, nonconvex penalties such as our RDP-based penalization scheme, may introduce local optima that are difficult to avoid if the initialization is chosen poorly.  Although in both these cases a better initialization typically yields better performance, it is undesirable to spend significant computational resources in doing so.

In practice, one approach is to initialize with an appropriately-scaled $A^T y$.  We have found this approach particularly effective in compressive sensing contexts where the sensing matrix $A$ acts as a near-isometry.  In many applications, a Fourier-based inversion scheme often leads to an effective initialization.  In the emission tomography example we consider in Sec.~\ref{sec:experiments}, a filtered back-projection estimate provides a sufficiently good initialization with low computational cost.  In particular, the initialization we use in the numerical experiments results from two iterations of a non-convergent version of the EPL-INC-3 algorithm, initialized by the filtered back-projection estimate. The EPL-INC-3 method employs an incremental penalized Poisson likelihood EM algorithm with an image roughness penalty, more details are available in Section~\ref{sec:algorithmsetup}.

\subsection{Termination criteria}
\label{sec:termination}

In this section, we list criteria by which \spiral\  decides whether the
iterates in the subproblem (\ref{eq:spiral_sub}) and for the main
problem (\ref{eq:spiral_goal}) are acceptable approximations to the true
minimizers to terminate the algorithm.  Here, we only provide criteria
for the \spiral-$\ell_1$ subproblem since the global solution to the
partition-based \spiral~subproblem (\ref{eq:partition}) can be easily
and exactly obtained using a non-iterative tree-pruning
algorithm \cite{willett:density}, even though its objective function is
nonconvex (due to the nonconvex penalty). For the TV-based method, 
we use the standard convergence criteria implemented by Beck and Teboulle \cite{BeckTV}.

\subsubsection{\spiral-$\ell_1$ subproblem}

The criterion for termination for the \spiral-$\ell_1$ subproblem measures the duality gap.  
Since the objective function $\phi^{k}(\theta)$
in (\ref{eq:sparsa_con}) is convex and all the constraints are affine, then (a weaker)
Slater's condition holds \cite{convexoptimization}  and, therefore, the duality gap is zero, i.e.,
$$
	\phi^{k}(\theta^{k+1}) = -h(\gamma^{\star}, \lambda^{\star}),
$$	
where $(\gamma^{\star}, \lambda^{\star})$ solves (\ref{eq:dual}).  Recall that, at the
$j\th$ iterate, $\theta^{(j)} = s^k + \gamma^{(j)} + W^T\lambda^{(j)}$ can be viewed as an approximate solution to $\theta^{k+1}$.  Thus, the \spiral\ approach for the $\ell_1$ subproblem  will consider the iterates to be sufficiently close to the optimal solution if
$$
	\frac{ | \phi^{k}(\theta^{(j)}) + h(\gamma^{(j)}, \lambda^{(j)}) |}{ | \phi^{k}(\theta^{(j)})|}
	\le \texttt{tol}_{\texttt{SUB}},
$$
where $\texttt{tol}_{\texttt{SUB}} > 0$ is some small constant.   In our numerical experiments, we often found that it is not necessary to
solve this subproblem very accurately, especially at the beginning of the algorithm where the iterate
$\theta^k$ is still far from the optimal solution.

\subsubsection{\spiral}

Since the global minimizer $f^{\star}$ of (\ref{eq:spiral_goal}) is not known
\emph{a priori}, criteria to terminate the \spiral\ algorithm must be established 
to determine whether a computed minimizer $\hf$ is an acceptable solution.
We list two such criteria.
The first of these criteria is simple: terminate if consecutive iterates or the corresponding
objective values do not change significantly, i.e., 
\begin{equation}
	{\| f^{k+1} - f^{k} \|_2}/{\| f^{k} \|_2} \le \texttt{tolP}
	\label{eq:spiraltermcriteria}
\end{equation}	
or
\begin{equation*}
	{| \Phi(f^{k+1}) - \Phi(f^k) |}/{| \Phi(f^k) |} \le \texttt{tolP},	
\end{equation*}
where \texttt{tolP} is a small positive constant.  The advantage of these criteria is that
they apply to general penalty functions.  The disadvantage, however, is that 
it is possible that the change between two consecutive iterates may be small or that they result in
only small improvements in the objective function even though iterates are still far from
the true solution.  However, we have yet to observe this premature termination in practice. 

The next criterion applies only to \spiral-$\ell_1$, where the penalty is convex and, after
a change of variables, differentiable.   This criterion is based on 
the Karush-Kuhn-Tucker  (KKT) conditions for optimality: at the $k\th$ iteration, given
$\theta^k$ and the corresponding Lagrange multipliers $\lambda^k$ computed 
from (\ref{eq:dual}), we determine whether
\begin{equation}\label{eq:crit1}
	\| \nabla \Phi(\theta^k) - W^T \lambda^k \|_2 
	\le \texttt{tolP},
\end{equation}
The left hand side corresponds to the gradient of the Lagrangian
function, which criterion (\ref{eq:crit1}) forces to be sufficiently
close to zero.  A complementarity condition could also be required,
but by construction, it is always satisfied, i.e., $(\lambda^k)^T
W\theta^k = 0$ using \eqref{eq:lambdaj} and \eqref{eq:Wthetaj}.

\subsection{Convergence proof}
\label{sec:convergence}

In this section we consider the minimization problem
\begin{equation}\begin{aligned}
	&\minimize_{f \in \Re^n} & &\Phi(f) \deq F(f) + \tau \pen(f) \\
	&\st & &f \ge 0,
\end{aligned}\end{equation}
where $F$ is the negative Poisson log-likelihood defined in \eqref{eq:negpoisson}. To ease analysis, we consider the equivalent problem of
\begin{equation}
	\minimize_{f \in \Re^n} F(f) + \rho(f),
	\label{eq:convergence_objective}
\end{equation}
where $\rho: \Re^n \to \overline{\Re} = \Re \cup \{-\infty,\infty\}$ is defined to be 
\begin{equation}
	\rho(f) = \tau \pen(f) + \delta_+(f),
\end{equation}
with $\delta_+ = \delta_{\Re_+^n}$ being the indicator function of the nonnegative orthant (our feasible set): 
\begin{equation}
\delta_+(f) = \begin{cases}0 & \text{ if $f \ge 0$}, \\ \infty & \text{ otherwise}.\end{cases}
\end{equation}
For notational simplicity, we still refer to the objective as $\Phi(f) = F(f) + \rho(f)$.
Considering the constraints in this manner allows one to avoid explicit examination of the convergence of the Karush-Kuhn-Tucker (KKT) conditions, and only consider subgradient calculus. For more details, see Rockafellar and Wets \cite{variationalanalysis}.

We make the following mild assumptions:
\begin{itemize}
\item (A1) $F$ is proper convex (i.e., $F$ is convex with $F(f) > -\infty$ for all $f$ and $F(f) < \infty$ for some $f$) and Lipschitz continuously differentiable on $\Re_+^n$,
\item (A2) $\rho$ is proper convex and continuous on $\Re_+^n$,
\item (A3) $\Phi$ is coercive (i.e., $\lim_{\|f\| \to \infty} \Phi(f) = \infty$),
\end{itemize}

We show that the assumption (A1) is satisfied by the negative Poisson log-likelihood through the following lemma.
\begin{lemma}
\label{lem:lipschitz}
The negative Poisson log-likelihood with parameter $\beta > 0$ is Lipschitz continuously differentiable on $\Re_+^n$ with Lipschitz constant 
\begin{equation*}
L \leq \frac{\max(y)}{\beta^2} \|A\|_2^2 \leq \frac{\max(y)}{\beta^2} \max(A^T \ones) \max(A\ones).
\end{equation*}
\end{lemma}
\begin{IEEEproof}
  Clearly the Hessian of $F$ \eqref{eq:hessian} is positive
  semidefinite; therefore all that is required is a bound on the
  largest eigenvalue of $\nabla^2 F(f)$ over the feasible set.  That
  is, we need to bound
\begin{equation*}
\lambda_{\rm max} = \sup_{f \geq 0} \sup_{\|z\|_2 \leq 1} z^T \nabla^2 F(f) z.
\end{equation*}
Since $A$ is nonnegative, the supremum over $f \geq 0$ is attained at $f = 0$, as this minimizes the denominator of the fraction in \eqref{eq:hessian}.  Therefore we simply need to bound the largest eigenvalue of  $\nabla^2 F(0) = \frac{1}{\beta^2} A^T \text{Diag}(y)A$. 
Using properties of matrix norms, we have
\begin{align*}
\lambda_{\rm max} &= \frac{1}{\beta^2} \| A^T \text{Diag}(y) A\|_2  
\leq \frac{1}{\beta^2} \|A\|_2^2 \|\text{Diag}(y)\|_2 \\
&= \frac{\max(y)}{\beta^2} \|A\|_2^2
\leq \frac{\max(y)}{\beta^2} \|A\|_1 \|A\|_\infty \\
&= \frac{\max(y)}{\beta^2} \max(A^T\ones) \max(A\ones).
\end{align*}
Since a bound on the largest eigenvalue of the Hessian is the same as a bound on the Lipschitz constant of the gradient, this completes the proof.
\end{IEEEproof}
We present the latter bound in Lemma \ref{lem:lipschitz} since the matrix $A$ is often such that computing $\|A\|_2$ may be difficult, yet matrix-vector products with $A$ and $A^T$ can be computed at far lower cost.

Before launching into the details of our convergence proof, let us recall a few facts.  A point $\overline{f}$ is said to be critical for \eqref{eq:convergence_objective} if 
\begin{equation*}
0 \in \partial\Phi(\overline{f}) = \nabla F(\overline{f}) + \partial \rho(\overline{f}),
\end{equation*}
where $\partial g(\overline{f})$ denotes the subdifferential (set of all subgradients) of $g$ at $\overline{f}$ \cite{variationalanalysis}. A minimizer $\overline{f}$ and corresponding minimum $\Phi(\overline{f})$ are guaranteed to exist for $\Phi$ coercive, proper convex, and lower semicontinuous, with criticality of $\overline{f}$ sufficient for $\overline{f}$ to optimize \eqref{eq:convergence_objective}. In solving \eqref{eq:convergence_objective}, we generate a sequence of iterates $\{f^k\}_{k \in \ints_+}$, which due to the nonnegativity constraint, we have $f^k \ge 0$ for all $k \in \ints_+$. 

We also simplify notation. Recall $\delta^{k+1} = f^{k+1} - f^{k}$, and let
$l(k) = \argmax_{i=[k-M]_+,\dots,k} \Phi(f^i)$,
i.e., the index where $\Phi$ is largest over the past $M$ iterations.  Then \eqref{eq:spiral_accept} is more simply written as
\begin{equation}
\Phi(f^{k+1}) \le \Phi(f^{l(k)}) - \tfrac{\sigma\alpha_k}{2}\|\delta^{k+1}\|_2^2.
\label{eq:convergence_acceptance}
\end{equation}
It can be shown that the sequence $\Phi(f^{l(k)})$ is nonincreasing, and further that $\Phi(f^0) = \Phi(f^{l(0)}) \ge \Phi(f^{l(k)}) \ge \Phi(f^{k})$.  Hence from assumption (A3), all the iterates are contained in a compact convex subset of $\Re_+^n$ (namely the initial sublevelset $\{f : \Phi(f) \le \Phi(f^0)\}$), and, therefore, we are guaranteed to have at least one convergent subsequence of iterates. 

Following the construction in \cite{sparsa}, the proof of our global convergence result (stated in Theorem \ref{thm:convergence} below) is best presented by first proving a set of supporting lemmas. The proofs of these lemmas are inspired by \cite{sparsa} and \cite{hager}, which consider the unconstrained setting only. Many important results hinge only on the analysis of the objective values and the acceptance rule \eqref{eq:convergence_acceptance}, and can be carried over to the constrained setting without any modification. In such cases we only highlight what restrictions must be in effect. 

However, some care must be taken in restricting the feasible domain to $\Re_+^n$ in the proof of Lemma~\ref{lem:stepifnotcritical}.
The proof in \cite{sparsa} relies on a result in \cite{convex}, which holds only for open sets, not the closed set $\Re_+^n$. To show this result in the constrained setting, we assume that $\rho$ is continuous on $\Re_+^n$. This assumption is stronger than the lower semicontinuity assumed in the unconstrained case; however it is still satisfied by the penalization schemes considered in Section~\ref{sec:algorithms}.

\begin{lemma}
\label{lem:stepifnotcritical}
Suppose $\overline{f} \in \Re_+^n$ is not critical for \eqref{eq:convergence_objective}. Then for any $\overline{\alpha} \ge \alpha_{\min}$, there exists an $\epsilon > 0$ such that for any subsequence $\{f^{k_j}\}_{j \in \ints_+}$ with $\lim_{j \to \infty} f^{k_j} = \overline{f}$, $f^{k_j} \in \Re_+^n$ for all $j \in \ints_+$, and $\alpha_{\min} \le \alpha_{k_j} \le \overline{\alpha}$, we have
\begin{equation*}
\| \delta^{k_j + 1} \|_2 = \|f^{k_j + 1} - f^{k_j}\|_2 \ge \epsilon
\end{equation*}
for all $j$ sufficiently large. 
\end{lemma}
\begin{IEEEproof}
Assume $\| \delta^{k_j + 1} \|_2 \to 0$ for contradiction, implying 
$\lim_{j\to\infty} f^{k_j + 1} = \lim_{j\to\infty} (\delta^{k_j + 1} + f^{k_j}) = \overline{f}$, where $f^{k_j + 1} \in \Re_+^n$ for all $j \in \ints_+$. Now since $f^{k_j + 1}$ is critical for the subproblem \eqref{eq:spiral_sub_general}, we have
\begin{equation*}
0 \in \nabla F(f^{k_j}) + \alpha_{k_j} \delta^{k_j + 1} + \partial \rho(f^{k_j+1}).
\end{equation*}
From the definition of $\partial \rho(f^{k_j+1})$, we have 
\begin{equation*}
\rho(z) \ge \rho(f^{k_j+1}) - [\grad F(f^{k_j}) + \alpha_{k_j} \delta^{k_j + 1}]^T(z - f^{k_j + 1}).
\end{equation*}
for all $z \in \Re^n$.
We now take the limit as $j \to \infty$.  Since both $\rho$ and $\grad F$ are continuous on $\Re_+^n$, $\rho(f^{k_j + 1}) \to \rho(\overline{f})$ and $\grad F(f^{k_j}) \to \grad F(\overline{f})$. Additionally, since $\delta^{k_j + 1}\to 0$ with $\alpha_{k_j}$ bounded, $\alpha_{k_j} \delta^{k_j + 1} \to 0$, hence for all $z \in \Re^n$ we have
\begin{equation*}
\rho(z) \ge \rho(\overline{f}) - \grad F(\overline{f})^T(z - \overline{f}),
\end{equation*}
showing $0 \in \grad F(\overline{f}) + \partial\rho(\overline{f})$, implying $\overline{f}$ is a critical point for \eqref{eq:convergence_objective}, the desired contradiction. Hence we cannot have $\|\delta^{k_j + 1}\|_2 \to 0$ as assumed. 
\end{IEEEproof}

\begin{lemma}
\label{lem:acceptancesatisfied}
Given $\sigma \in (0,1)$, then there exists an $\widetilde{\alpha} > 0$ such that for any sequence $\{f^i\}_{i=1}^k$ with $f^i \in \Re_+^n$ for all $i = \ints_+$, the acceptance criterion \eqref{eq:convergence_acceptance} is satisfied whenever $\alpha_k \ge \widetilde{\alpha} = 2L/(1-\sigma)$, where $L$ is the Lipschitz constant for $\grad F$ over $\Re_+^n$.
\end{lemma}
\begin{IEEEproof}
The proof follows similarly from the proof of Lemma 4 in \cite{sparsa}, where we only consider $f^i \in \Re_+^n$, and hence only require $F$ to be Lipschitz continuously differentiable over $\Re_+^n$.
\end{IEEEproof}

\begin{lemma}
\label{lem:steptozero}
The sequence $\{f^{k}\}_{k \in \ints_+}$ generated by SPIRAL is such that $\lim_{k\to\infty} \delta^{k+1} = 0$, and there exists $\overline{\Phi} \in \Re$ such that $\lim_{k\to\infty} \Phi(f^k) = \overline{\Phi}$.
\end{lemma}

\begin{IEEEproof} The proof follows identically as the proof of Lemma 4 in \cite{sparsa} as no modifications are required for the constraints within the proof.
\end{IEEEproof} 

With the above lemmas in hand, we are now ready to state the convergence result.
\begin{theorem}
\label{thm:convergence}
Suppose that SPIRAL is applied to \eqref{eq:convergence_objective}, where (A1) through (A3) hold, then all accumulation points are critical points, and hence SPIRAL converges to a minimizer of \eqref{eq:convergence_objective}. Moreover, the sequence of objective values converges sublinearly to the minimal value.  Additionally, if $F$ is strongly convex, the sequence of objective values converges $R$-linearly. 
\end{theorem}

\begin{IEEEproof}
Using the above lemmas, the proof that the objective values converge follows from \cite{sparsa}. 
Assume for contradiction that $\overline{f}$ is an accumulation point that is not critical. Let $\{f^{k_j}\}_{j \in \ints_+}$ a subsequence such that $\lim_{j\to\infty}f^{k_j} = \overline{f}$. If the sequence $\{\alpha_{k_j}\}_{j \in \ints_+}$ were bounded, we would have from Lemma \ref{lem:stepifnotcritical} that $\|\delta^{k_j+1}\|_2 \ge \epsilon$ for some $\epsilon > 0$ and all $j$ large enough. However, this is in contradiction Lemma \ref{lem:steptozero}, hence $\{\alpha_{k_j}\}_{j \in \ints_+}$ must be unbounded. In particular, we must have that for some $j$ large enough, $\alpha_{k_j} \ge \eta \max(\alpha_{\max}, \widetilde{\alpha})$, meaning $\alpha = \alpha_{k_j}/\eta \ge \widetilde{\alpha}$ must have been tried yet failed the acceptance criterion, which is prohibited by Lemma \ref{lem:acceptancesatisfied}, since any $\alpha \ge \widetilde{\alpha}$ satisfies the acceptance criterion. This further contradiction shows that any noncritical point cannot be an accumulation point, hence all critical points are accumulation points, and due to the convexity of $\Phi$, $\overline{f}$ optimizes \eqref{eq:convergence_objective}.

Following \cite{hager}, it can be shown that for some constant $c \ge 0$
\begin{equation*}
\Phi(f^{k}) - \Phi(\overline{f}) \le \frac{c}{k},
\end{equation*}
that is the objective values converge at a sublinear rate.  Moreover, if $F$ is strongly convex, then there exist constants $C > 0$ and $r \in (0,1)$ such that
\begin{equation*}
\Phi(f^{k}) - \Phi(\overline{f}) \le Cr^{k}\left(\Phi(f^0) - \Phi(\overline{f})\right),
\end{equation*}
hence the objective values exhibit $R$-linear convergence.
\end{IEEEproof}

\subsection{Uniqueness of the minimizer}
\label{sec:uniqueness}

Under the general assumptions (A1)--(A3) in Section~\ref{sec:convergence}, the set $\argmin_f \Phi(f)$ will be a compact convex set (c.f. \cite{variationalanalysis}) and may not be the singleton $\{\hf\}$. That is, the solution to the minimization \eqref{eq:spiral_goal} may not be unique.
The most general condition under which the minimizer $\hf$ is unique is when the objective $\Phi$ is strictly convex, which is obtained when either $\phi$ or the penalty $\pen$ is strictly convex. In the underdetermined case of interest, $\phi$ cannot be strictly convex, since it is trivial to find a vector $z$ in the kernel of $A$, and as a consequence $F(f + \alpha z) = F(f)$ for any $\alpha \in \Re$. In this case, a strictly convex objective can be obtained if the penalty is strictly convex, however most interesting choices for the penalty term (e.g., $\ell_1$-norm, total variation) preclude this possibility.  More precisely, the penalty need only be strictly convex on the null space of $A$, however this condition is difficult to verify in practice.

A different approach must then be taken which utilizes the precise choice of $\pen$ and functional form for $F$, in which conditions for a unique minimizer are established on a case-by-case basis. We do not consider the RDP-based penalty, since convergence to a global optimum cannot be explicitly guaranteed, and therefore only focus on the $\ell_1$ and TV penalties. We show for these two penalties that if the solution is not unique, then there exist two distinct solutions $\widehat{f}$ and $\widetilde{f}$ such that $A\widehat{f} = A\widetilde{f}$ and $\pen(\widehat{f}) = \pen(\widetilde{f})$. This means that these solutions are identical in the sense that they are equally faithful to the data and equally sparse or smooth as measured by the penalty. In certain cases, $A$ is such that $A\widehat{f} = A\widetilde{f}$ only if $\widehat{f} = \widetilde{f}$, hence the solution is unique. 

We begin the analysis by expressing the TV seminorm as $\|f\|_{\text{TV}} = \|Df\|_1$, with $D = [D_1 ; D_2]$
where $D_1$ and $D_2$ are respectively the horizontal and vertical first-order difference matrices.  Therefore both the $\ell_1$ and TV penalties can be cast into the form $\pen(f) = \|Bf\|_1$ where $B = \tau W^T$ for the $\ell_1$ penalty, and $B = \tau D$ for the TV penalty. The KKT optimality conditions for a solution $\widehat{f}$ are such that there exists a corresponding Lagrange multiplier vector $\widehat{\lambda}$ (not necessarily unique) such that, together with $\widehat{f}$,
\begin{equation}\begin{aligned}
0 &\in \grad F(\widehat{f}) + \partial \pen(\widehat{f}) - \widehat{\lambda}, & \widehat{f} &\ge 0, \\
0 &= \widehat{\lambda}^T \widehat{f}, & \widehat{\lambda} &\ge 0.
\label{eq:uniquekkt}
\end{aligned}\end{equation}
To examine the subdifferential of the penalty, note the composition rule that if $g(x) = h(Bx)$, then $\partial g(x) = B^T \partial h(Bx)$, that is if $b_i$ are the rows of $B$, then
\begin{equation*}
\partial g(x) = \Big\{\sum_{i} s_i b_i, s_i \in (\partial h(Bx))_i\Big\}.
\end{equation*}
When $h(x) = \|x\|_1$, the subdifferential is
\begin{equation*}
(\partial h(x))_i = \begin{cases} \text{sign}(x_i) & \text{ if $x_i \ne 0$}, \\ [-1, 1] & \text{ if $x_i = 0$},\end{cases}
\end{equation*}
and thus the first KKT condition in \eqref{eq:uniquekkt} is equivalent to
\begin{equation*}
0 \in \Big\{ \grad F(\widehat{f}) - \widehat{\lambda} + 
\sum_{i \in \widehat{S}} \text{sign}(b_i^T \widehat{f}) b_i + \sum_{i \notin \widehat{S}} s_i b_i : s_i \in [-1, 1] \Big\}
\end{equation*}
where $\widehat{S} = \{ i : b_i^T \widehat{f} \ne 0 \}$ is the support set of $B\widehat{f}$.

If $\widehat{f}$ is not the unique solution, then there exists another solution $\widetilde{f}$, distinct from $\widehat{f}$, that also satisfies the KKT conditions with a corresponding $\widetilde{\lambda}$.  What we show next is that, without loss of generality, we can assume that $\widehat{f}$ and $\widetilde{f}$ share certain properties that simplifies the proceeding analysis. 

We know that $\widehat{f}$ and $\widetilde{f}$ must lie on a compact convex set. 
In particular, any point
$z(\gamma) = (1-\gamma) \widehat{f} + \gamma \widetilde{f}$ with $\gamma \in [0, 1]$
must also be a solution. Define the following index sets:
\begin{equation*}\begin{aligned}
S(x) &= \{i : b_i^T x \ne 0\}, & Q(x) &= \{i : x_i \ne 0\}.
\end{aligned}\end{equation*}
Now clearly $S(z(\gamma)) \subseteq S(\widehat{f}) \cup S(\widetilde{f})$ and $Q(z(\gamma)) \subseteq Q(\widehat{f}) \cup Q(\widetilde{f})$, but as we sweep through $\gamma$, the sets $S(z(\gamma))$ and $Q(z(\gamma))$ may not remain constant.  However they will be constant on certain intervals. To see this, define the following:
\begin{equation*}\begin{aligned}
\Gamma_S &= \{\gamma \in [0,1] : b_i^Tz(\gamma) = 0 \text{ for some } i \in S(\widehat{f}) \cup S(\widetilde{f})\},\\
\Gamma_Q &= \{\gamma \in [0,1] : z_i(\gamma) = 0 \text{ for some } i \in Q(\widehat{f}) \cup Q(\widetilde{f})\},\\
\Gamma &= \Gamma_S \cup \Gamma_Q.
\end{aligned}\end{equation*}
These sets basically count the number of zero-crossings that occur in $z(\gamma)$, excluding the components that are always zero (e.g,, if the line segment between $\widehat{f}$ and $\widetilde{f}$ lies on a coordinate plane), and hence they are finite sets. Therefore, $\Gamma = \{\gamma_i\}_{i=1}^N$ for some $N < \infty$, hence there exist open intervals $(\gamma_i,\gamma_{i+1})$ on which $S(z(\gamma))$ and $Q(z(\gamma))$ are constant. So without loss of generality, we can select $\widehat{f}$ and $\widetilde{f}$ to be two distinct points in one of these intervals, so that we now have
 $S \triangleq S(\widehat{f}) = S(\widetilde{f})$ and $Q \triangleq Q(\widehat{f}) = Q(\widetilde{f})$. Furthermore, we will also have $\text{sign}(b_i^T \widehat{f}) = \text{sign}(b_i^T \widetilde{f})$ for all $i \in S$.
 
Using the properties established above, and from the first KKT condition for both $(\widehat{f},\widehat{\lambda})$ and $(\widetilde{f},\widetilde{\lambda})$, we have that
\begin{equation*}\begin{aligned}
0 \in \Big\{ &\grad F(\widetilde{f}) - \grad F(\widehat{f}) - \widetilde{\lambda} + \widehat{\lambda} \\
&+ \sum_{i \notin S}(\widetilde{s}_i - \widehat{s}_i)b_i : \widetilde{s}_i, \widehat{s}_i \in [-1, 1] \Big\}.
\end{aligned}\end{equation*}
From here we multiply componentwise by $(\widetilde{f} - \widehat{f})$ and sum to yield
\begin{equation*}
0 = (\widetilde{f} - \widehat{f})^T(\grad F(\widetilde{f}) - \grad F(\widehat{f})) + 
\widehat{f}^T\widetilde{\lambda} + \widetilde{f}^T\widehat{\lambda},
\end{equation*}
where we have used the complementarity KKT condition that $\widehat{f}^T\widehat{\lambda} = \widetilde{f}^T\widetilde{\lambda} = 0$, and the fact that $b_i^T\widehat{f} = b_i^T\widetilde{f} = 0$ for all $i \notin S$. Now since $\widetilde{\lambda}_i > 0$ only if $\widetilde{f}_i = 0$ (similarly $\widehat{\lambda}_i > 0$ only if $\widetilde{f}_i = 0$), and since $\widetilde{f}_i = \widehat{f}_i = 0$ for all $i \notin Q$, we also have that $\widehat{f}^T\widetilde{\lambda} = \widetilde{f}^T\widehat{\lambda} = 0$. Therefore we deduce that 
\begin{equation*}
0 = \sum_{i=1}^m \frac{y_i}{(e_i^TA\widetilde{f} + \beta)(e_i^TA\widehat{f} + \beta)} \left[ e_i^TA(\widetilde{f} - \widehat{f})\right]^2,
\end{equation*}
and since this is a sum of nonnegative terms (if $y > 0$, the fraction is strictly positive) it can only assume a zero value if $A(\widehat{f} - \widetilde{f}) = 0$. This means that if both $\widehat{f}$ and $\widetilde{f}$ are solutions, they need to be equivalent up to their projections. Since they are both minimizers of the objective, they necessarily have identical objective values, and since they both achieve the same value for $F$, it is clear that $\pen(\widehat{f}) = \pen(\widetilde{f})$. Hence the solutions are also equally sparse or smooth with respect to the penalty.

In the case of the $\ell_1$ norm, we can be a bit more precise by considering the fact that we often examine the solution in terms of the coefficients $\widehat{\theta} = W^T \widehat{f}$. In this case, if $AW$ satisfies a restricted-isometry type property \cite{RIP} for vectors of sparsity $s = 2|S|$, we have $(1-\delta_s)\|\widehat{\theta} - \widetilde{\theta}\|_2 \le \|AW(\widehat{\theta} - \widetilde{\theta})\|_2$, and hence $AW(\widehat{\theta} - \widetilde{\theta}) = 0$ only if $\widehat{\theta} - \widetilde{\theta} = 0$, contradicting the distinctness of $\widehat{\theta}$ and $\widetilde{\theta}$, hence if $AW$ satisfies RIP, the solution is unique.

\section{Numerical experiments}
\label{sec:experiments}

\subsection{Simulation Setup}

Although the algorithms described heretofore are applicable to a wide range of imaging contexts, here we demonstrate the effectiveness of the proposed methods on a simulated limited-angle emission tomography dataset.  We compare our algorithm with the currently available state-of-the-art emission tomography reconstruction algorithms \cite{fesslerSPSOS}. 

\begin{figure}[t]
	\centering
	\subfigure[]{\includegraphics[width=0.83in]{GroundTruth.png}}
	\subfigure[]{\includegraphics[width=0.83in]{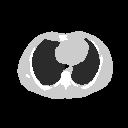}}
	\subfigure[]{\includegraphics[width=0.83in]{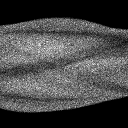}}
	\subfigure[]{\includegraphics[width=0.83in]{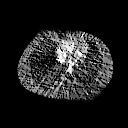}}
	\caption{Experimental setup:  (a) true emission image, (b) attenuation map,
		(c) noisy projection data, (d) estimate used as initialization.}
	\label{fig:ExpSetup}
\end{figure}

In this experimental setup, we wish to reconstruct the true axial emission map $\ftrue$ (Fig.~\ref{fig:ExpSetup}(a)) as accurately as possible. The photon flux described by this emission map is subject to the attenuation effects caused by the various densities of tissue through which the photons must travel to reach the detector array.  The simulated attenuation map $\mu$ (Fig.~\ref{fig:ExpSetup}(b)) is assumed known during the reconstruction process.  The simulated emission and attenuation images are standard test images included in the Image Reconstruction Toolbox (IRT) by Fessler \cite{fesslerIRT}.

The limited-angle tomographic projection $R$, corresponding to
parallel strip-integral geometry with 128 radial samples and 128
angular samples spaced uniformly over 135 degrees, was also generated
by the IRT software \cite{fesslerIRT}.  The resulting sensing matrix
is then given by $A \equiv {\rm diag}[\exp(-R\mu)]R$, with the noisy
tomographic data $y$ simulated according to the inhomogeneous Poisson
process \eqref{eq:PoissonModel}.  We simulated ten realizations of the
data $y$ in order to examine the ten-trial average performance of all
the reconstruction methods presented.  We only show images
reconstructed using a particular realization of the data shown in
Fig.~\ref{fig:ExpSetup}(c)).  In this case, the noisy sinogram
observations have a total photon count of ${\rm 2.0\times10^5}$, a
mean count over the support of the tomographic projections of 18.08,
and a maximum count of 44.

\subsection{Algorithm Setup} 
\label{sec:algorithmsetup}

We evaluate the proposed \spiral\ approaches ($\ell_1$, TV, RDP, and
RDP-TI), two competing Poisson reconstruction methods
\cite{fesslerSPSOS,fesslerIRT}, and the SpaRSA algorithm
\cite{sparsa}. For the \spiral-$\ell_1$ method, the Daubechies family
(DB-2 through DB-8) were evaluated as the basis $W$. Here we only
present results using the DB-6 wavelet basis as the results were
similar across bases with DB-6 having a marginal lead in
performance. Using the \spiral-TV method, we investigate the effects
of not enforcing the near-monotonicity acceptance criterion
\eqref{eq:spiral_accept}.  For all methods that enforce
\eqref{eq:spiral_accept}, we used $\eta = 2$, $\sigma = 0.1$, $M =
10$.  Also, as described in Section~\ref{sec:termination}, we evaluate
how the accuracy of solving the subproblem impacts the global
performance of the \spiral-$\ell_1$ and TV methods.  When solving
these subproblems stringently (denoted Tight in the figures and
tables), a minimum of 10 and a maximum of 100 iterations were used,
with a convergence tolerance $\texttt{tol}_{\texttt{SUB}} = {\rm
  1\times10^{-8}}$; these parameters were relaxed to a maximum of 10
iterations, and a tolerance $\texttt{tol}_{\texttt{SUB}} = {\rm
  1\times10^{-4}}$ when the subproblems are solved less exactly
(denoted Loose). In all the \spiral\ methods, we set
$\alpha_\text{max} = 1/\alpha_\text{min} = {\rm 1 \times 10^{30}}$.

\newlength{\reconfigsize}
\setlength{\reconfigsize}{0.975in}
\begin{figure*}[t]
	\centering
	\subfigure[Ground Truth]
		{\includegraphics[width=\reconfigsize]{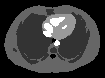}}
	\subfigure[\spiral-TV \hspace{\reconfigsize} Loose Monotonic
		\hspace{\reconfigsize} RMSE = 24.404\%]
		{\includegraphics[width=\reconfigsize]{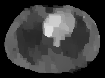}}
	\subfigure[\spiral-TV \hspace{\reconfigsize} Loose Nonmonotonic
		\hspace{\reconfigsize} RMSE = 24.962\%]
		{\includegraphics[width=\reconfigsize]{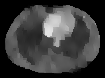}} 
	\subfigure[\spiral-TV \hspace{\reconfigsize} Tight Monotonic
		\hspace{\reconfigsize} RMSE = 24.526\%]
		{\includegraphics[width=\reconfigsize]{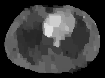}}
	\subfigure[\spiral-TV \hspace{\reconfigsize} Tight Nonmonotonic
		\hspace{\reconfigsize} RMSE = 24.467\%]
		{\includegraphics[width=\reconfigsize]{SPIRALTVLooseNonmono_crop.png}}		
	\subfigure[\spiral-RDP \hspace{\reconfigsize} (Translation Variant)
		\hspace{\reconfigsize} RMSE = 33.959\%]
		{\includegraphics[width=\reconfigsize]{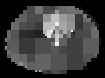}}
	\subfigure[\spiral-RDP-TI \hspace{\reconfigsize} (Cycle-Spun)
		\hspace{\reconfigsize} RMSE = 27.557\%]
		{\includegraphics[width=\reconfigsize]{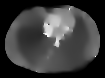}}		\\
	\subfigure[\spiral-$\ell_1$ \hspace{\reconfigsize} Loose (DB-6)
		\hspace{\reconfigsize} RMSE = 28.626\%]
		{\includegraphics[width=\reconfigsize]{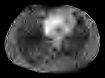}}
	\subfigure[\spiral-$\ell_1$ \hspace{\reconfigsize} Tight (DB-6)
		\hspace{\reconfigsize} RMSE = 28.665\%]
		{\includegraphics[width=\reconfigsize]{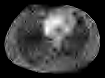}} 
	\subfigure[SpaRSA \hspace{\reconfigsize} $\ell_2$-$\ell_1$ (DB-6)
		\hspace{\reconfigsize} RMSE = 31.172\%]
		{\includegraphics[width=\reconfigsize]{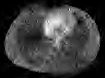}} 	
	\subfigure[SPS-OS \hspace{\reconfigsize} (Huber potential)
		\hspace{\reconfigsize} RMSE = 27.555\%]
		{\includegraphics[width=\reconfigsize]{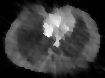}}
	\subfigure[SPS-OS \hspace{\reconfigsize} (Quadratic Potential)
		\hspace{\reconfigsize} RMSE = 29.420\%]
		{\includegraphics[width=\reconfigsize]{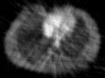}} 
	\subfigure[EPL-INC-3 \hspace{\reconfigsize} (Huber Potential)
		\hspace{\reconfigsize} RMSE = 24.748\%]
		{\includegraphics[width=\reconfigsize]{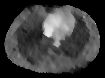}}
	\subfigure[EPL-INC-3 \hspace{\reconfigsize} (Quadratic Potential)
		\hspace{\reconfigsize} RMSE = 26.474\%]
		{\includegraphics[width=\reconfigsize]{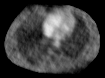}} 
	\caption{Single-trial reconstructed images for all methods considered.  
	Note $\text{{\rm RMSE (\%)}} = 100\cdot\|\hf - \ftrue\|_2/\|\ftrue\|_2$.}
	\label{fig:Reconstructions}
\end{figure*}

We compare our proposed approaches with two competing Poisson
reconstruction methods. The first, denoted SPS-OS, uses a separable
paraboloidal surrogate with ordered subsets algorithm
\cite{fesslerSPSOS}.  The second, denoted EPL-INC-3, employs an
incremental penalized Poisson likelihood EM algorithm and was
suggested by Prof.~Fessler as representative of the current
state-of-the-art in emission tomographic reconstruction.  Both of
these methods are available as part of the IRT \cite{fesslerIRT};
specifically, we used the \texttt{pwls\_sps\_os} and \texttt{epl\_inc}
functions from the toolbox.  In addition to these Poisson methods, we
also compare to the SpaRSA algorithm \cite{sparsa} which solves the
$\ell_1$-regularized least-squares ($\ell_2$-$\ell_1$) problem. Like
the \spiral-$\ell_1$ approaches, we only present results when $W$ is
the DB-6 wavelet basis. As the solution provided by SpaRSA is not
guaranteed to be nonnegative, we threshold the result to obtain a
feasible -- and therefore more accurate -- solution.  Including this
result allows us to demonstrate the effectiveness of solving the
formulation \eqref{eq:spiral_goal} that utilizes the Poisson
likelihood.

All of the methods considered here were initialized with the estimate
shown in Fig.~\ref{fig:ExpSetup}(d). This initialization results from
two iterations of a non-convergent version of the EPL-INC-3 algorithm,
itself initialized by the filtered back-projection estimate.  All
algorithms executed for a minimum of 50 iterations, and global
convergence was declared when the relative change in the iterates
\eqref{eq:spiraltermcriteria} fell below $\texttt{tolP} = {\rm 5
  \times 10^{-4}}$.

Lastly, in all of the experiments presented in this paper, we chose any parameters associated with each algorithm (such as $\tau$) to minimize the RMS error ($\text{{\rm RMSE (\%)}} = 100\cdot\|\hf - \ftrue\|_2/\|\ftrue\|_2$) of the reconstruction. While this would not be possible in practice, it does allow us to compare the {\em best-case} performance of various algorithms and penalization methods. In practical settings, regularization parameters can be chosen via cross-validation. This is particularly well-suited to many photon-limited imaging applications in which each detected photon has a time stamp associated with it; this timing information can be used to construct multiple independent and identically distributed realizations of the underlying Poisson process in software. The details of this procedure are a significant component of our ongoing research.

\subsection{Results analysis}

From the results presented in Table~\ref{tab:RMSEandTIME}, we see that
out of the proposed \spiral\ approaches, the method utilizing the
total variation penalization achieves the lowest RMSE and highest
visual quality, followed by \spiral-RDP-TI and \spiral-$\ell_1$.  This
bolsters the notion that there is much to be gained by considering
additional image structure beyond that of a sparse representation in
the basis $W$.  The \spiral-RDP method based on the non-cycle-spun
partitions simply are not competitive due to the high bias in
considering only a single shift of the RDP structure. Examining
Fig.~\ref{fig:Reconstructions}(f), we see that this bias is manifest
in the image as blocking artifacts due to the RDP structure not
fortuitously aligning to any strong edges in the image. Also from
Table~\ref{tab:RMSEandTIME} we see that the EPL-INC-3 method with the
Huber potential offers the toughest competition to the RMSE achieved
by the \spiral-TV approaches.

\begin{table}[b!]
	\centering
	\begin{tabular}{@{}lrrrr@{}}
		\toprule
										& \multicolumn{2}{c}{Single-Trial} 	& \multicolumn{2}{c}{Ten-Trial Average} 	\\
										\cmidrule(rl){2-3} 						\cmidrule(l){4-5}
		Method 						& RMSE (\%) 	& Time (s) 				& RMSE (\%) 	& Time (s) 			\\
		\midrule
		SPS-OS Huber 					& 27.555		&  9.078 					& 27.057 		&  9.094 				\\
		SPS-OS Quad 					& 29.420 		&  5.354 					& 29.049 		&  5.574 				\\
		EPL-INC-3 Huber 				& 24.748 		& 20.302 					& 24.462 		& 16.704 				\\
		EPL-INC-3 Quad 				& 26.474 		&  8.713 					& 26.013 		&  9.647 				\\
		SpaRSA 						& 31.172 		&  5.946 					& 29.987 		&  3.730 				\\
		\spiral-$\ell_1$ (L) 	 	& 28.626 		&  6.933 					& 28.050 		&  5.396 				\\
		\spiral-$\ell_1$ (T) 	 	& 28.665 		& 15.547 					& 27.980 		& 20.006 				\\
		\spiral-TV (L, M) 			& 24.404 		& 15.418 					& 24.270 		& 10.102 				\\
		\spiral-TV (L, NM) 			& 24.962 		&  7.821 					& 24.868 		&  4.721 				\\
		\spiral-TV (T, M) 			& 24.526 		& 25.423 					& 24.352 		& 20.900 				\\
		\spiral-TV (T, NM) 			& 24.467 		& 21.505 					& 24.571 		& 15.046 				\\
		\spiral-RDP 					& 33.959 		&  3.118 					& 34.946 		&  2.323 				\\
		\spiral-RDP-TI 				& 27.557 		&  5.313 					& 27.669 		&  4.756 				\\
		\bottomrule
	\end{tabular}
	\caption{Reconstruction RMSE and computation time for both the single-trial results, and results averaged over ten trials.
		Note: L = Loose, T = Tight, M = Monotonic, NM = Nonmonotonic, 
		$\text{{\rm RMSE (\%)}} = 100\cdot\|\hf - \ftrue\|_2/\|\ftrue\|_2$.}
	\label{tab:RMSEandTIME}
\end{table}

\begin{figure}[t!]
	\centering
	\includegraphics[scale=0.85]{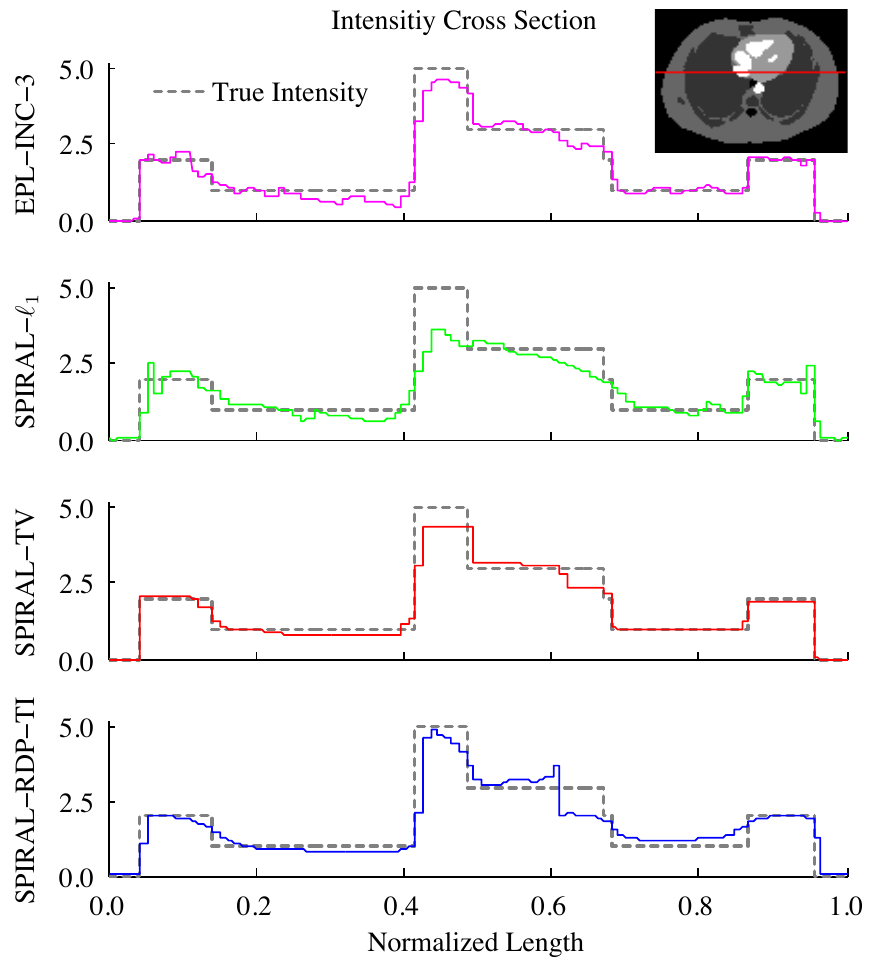}
	\caption{Intensity cross sections of the best-in-class methods: EPL-INC-3 with Huber potential, 
		\spiral-$\ell_1$ with loose subproblem convergence criteria, 
		monotonic \spiral-TV with loose subproblem convergence criteria,
		and \spiral-RDP-TI. The insert image shows the location of the cross section.  
		Note that the pixel boundaries in the image cause the stepwise nature of the profiles.}
	\label{fig:CrossSection}
	\vspace{-2ex}
\end{figure}

Examining the reconstructed images in Fig.~\ref{fig:Reconstructions},
we see that the visual fidelity correlates strongly with the RMSE.
The \spiral-TV and EPL-INC-3 results clearly do best at capturing
strong edges in the image while smoothing noise in the homogeneous
regions, with the \spiral-TV method more correctly smoothing
homogeneous regions.  This conclusion is supported by
Fig.~\ref{fig:CrossSection} in which we further analyze the
best-in-class \spiral\ and EPL-INC-3 methods. Between the EPL-INC-3
methods, the Huber potential is more effective at reducing the rough
noise-like deterioration that is quite pronounced when using the
quadratic potential. Here we see that both \spiral-TV and EPL-INC-3
results have sharp edges, but significant variations in the regions of
homogeneous intensity are a consistent source of error for the
EPL-INC-3 method. These variations cause many darker spots that could
be misinterpreted as regions of low uptake.
Figures~\ref{fig:Reconstructions} and \ref{fig:CrossSection} both show
that the \spiral-RDP-TI and $\ell_1$-based approaches tend to
over-smooth the entire image; however the \spiral-RDP-TI method better
captures the high intensity regions.  High-scale wavelet artifacts are
seen throughout the \spiral-$\ell_1$ and SpaRSA reconstructions. In
the SpaRSA results shown in this section, we clipped the final
reconstruction to be nonnegative. Even after this clipping operation,
the RMSE of the final result was significantly higher than the RMSE
associated with methods using the Poisson log-likelihood. This
suggests that the effort required to use the Poisson log-likelihood is
justified by tangible performance gains over methods associated with a
penalized least-squares objective.  Lastly, prominent streaking
artifacts in the SPS-OS methods hinder their accuracy, with using a
Huber potential resulting in a slightly better performance than a
quadratic potential.

\begin{figure}[t!]
	\centering
	\includegraphics[scale=0.85]{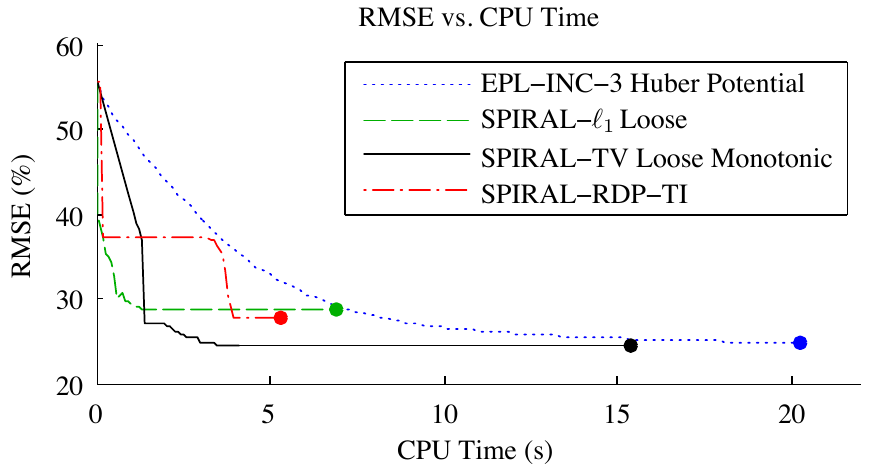}
	\caption{Convergence of the RMSE vs computation time for the best-in-class methods. The termination criteria
		was reached at the point indicated by the solid circle.}
	\label{fig:RMSEvsCPUbest}
	\vspace{-2ex}
\end{figure}

While accuracy and visual quality are held paramount in tomographic reconstruction, the computational cost and convergence behavior of the chosen methods should not be neglected.  Figure~\ref{fig:RMSEvsCPUbest}
shows the convergence behavior of the best-in-class methods.  Although the \spiral-$\ell_1$ and RDP-TI approaches yield higher RMSE, they show desirable convergence behavior with sharp initial decreases in RMSE and termination in under eight seconds. The \spiral-TV method also shows this sharp decrease, but exhibits a prolonged period where the RMSE changes little.  The EPL-INC-3 method shows none of these traits, slowly changing the RMSE throughout its execution until it finally reaches the termination criteria after 20 seconds.  

\begin{figure}[b!]
	\centering
	\includegraphics[scale=0.85]{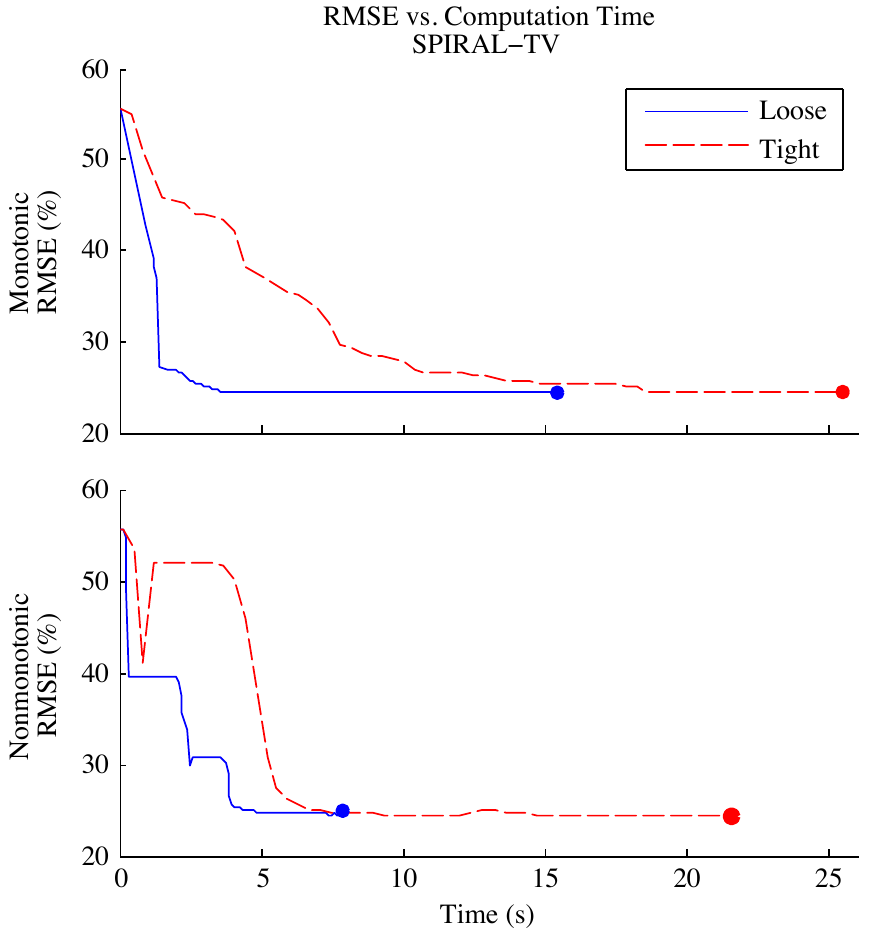}
	\caption{Comparison of the RMSE convergence rates for the different variations of the \spiral-TV method.
		The termination criteria was reached at the point indicated by the solid circle.}
	\label{fig:RMSEvsCPUtv}
	\vspace{-2ex}
\end{figure}

Next consider Fig.~\ref{fig:RMSEvsCPUtv}, where we show the RMSE convergence for the variations of the \spiral-TV algorithm.  Comparing the top and bottom axes, we see that by not enforcing the near-monotonicity acceptance criteria in \eqref{eq:spiral_accept} accelerates the convergence rate at the potential cost of large increases of the RMSE during the execution of the algorithm (seen at the near the two second mark when tightly solving the TV subproblem). Although the nonmonotonic algorithms typically yield early sharp decreases in RMSE, the ten-trial average behavior in Table~\ref{tab:RMSEandTIME} shows that the monotonic algorithm will achieve a more accurate solution at termination. This suggests that a two-stage approach of starting with a nonmonotonic algorithm to achieve quick, yet approximate, convergence then switching to the monotonic algorithm may be the best overall approach.

\begin{figure}[t!]
	\centering
	\includegraphics[scale=0.85]{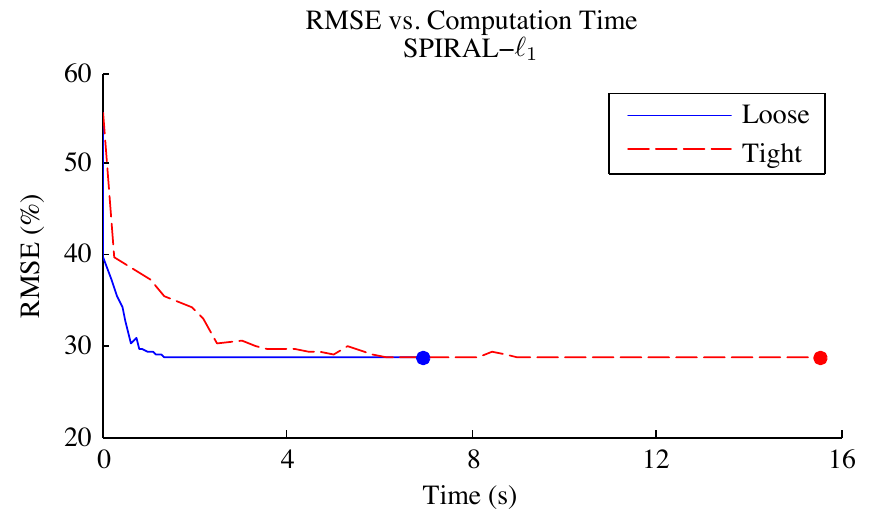}
	\caption{Comparison of the RMSE convergence rates for the different variations of the \spiral-$\ell_1$ method.
		The termination criteria was reached at the point indicated by the solid circle.}
	\label{fig:RMSEvsCPUl1}
	\vspace{-2ex}
\end{figure}

Figures~\ref{fig:RMSEvsCPUtv} and \ref{fig:RMSEvsCPUl1} support the
use of relaxed convergence criteria for the TV and $\ell_1$
subproblems, Table~\ref{tab:RMSEandTIME} also shows that little is
lost in terms of final RMSE by using this less-stringent criterion.
In fact, by combining the nonmonotonic option and loose subproblem
convergence criteria results in a \spiral-TV approach that, on
average, nearly matches the best EPL-INC-3 method in terms of RMSE,
yet is approximately four times faster.  Further, this algorithm is
even faster than the \emph{least} accurate SPS-OS method.  Therefore,
we believe that further investigation of nonmonotonic and approximate
subproblem algorithms is a fruitful area of research for fast and
accurate image reconstruction.  We finish the discussion of the
numerical results by directing the interested reader to additional
experimental results in \cite{PoissonCS_tsp,PoissonStreams} which
demonstrate the effectiveness of the proposed \spiral\ method in a
compressed sensing context.

\section{Conclusion}

\label{sec:conclusion}

We have formulated the general goal of reconstructing an image from
photon-limited measurements as a penalized maximum Poisson likelihood
estimation problem.  To obtain a solution to this problem, we have
proposed an algorithm that allows for a flexible choice of
penalization methods, and focused particularly on sparsity-promoting
penalties.  In particular, we detail the cases where the penalty
corresponds to the sparsity-promoting $\ell_1$ norm of the expansion
coefficients in a sparsifying basis, is related to the complexity of a
partition-based estimate, or is proportional to the total variation of
the image.  We establish mild conditions for which this algorithm has
desirable convergence properties, although in practice it is
beneficial to relax these conditions to attain faster (albeit
nonmonotonic) convergence.  We demonstrate the effectiveness of our
methods through a simulated emission tomography example.  When our
total variation regularized method is applied to this problem, the
resulting estimates outperform the current state-of-the-art approaches
developed specifically for emission tomography.  In particular, it
results in fewer spurious artifacts than wavelet-regularized methods,
and unlike partition-regularized methods, is the result of a convex
optimization procedure.

\section*{Acknowledgment}

This work was supported by NSF CAREER Award No.\
    CCF-06-43947, NSF Award No.\ DMS-08-11062, DARPA Grant
    No. HR0011-07-1-003, and AFRL Grant No. FA8650-07-D-1221.

\ifCLASSOPTIONcaptionsoff
\fi



\bibliographystyle{IEEEtran}
\bibliography{PoissonCS}
%



%

\begin{IEEEbiography}[{\includegraphics[width=1in,height=1.25in,clip,keepaspectratio]{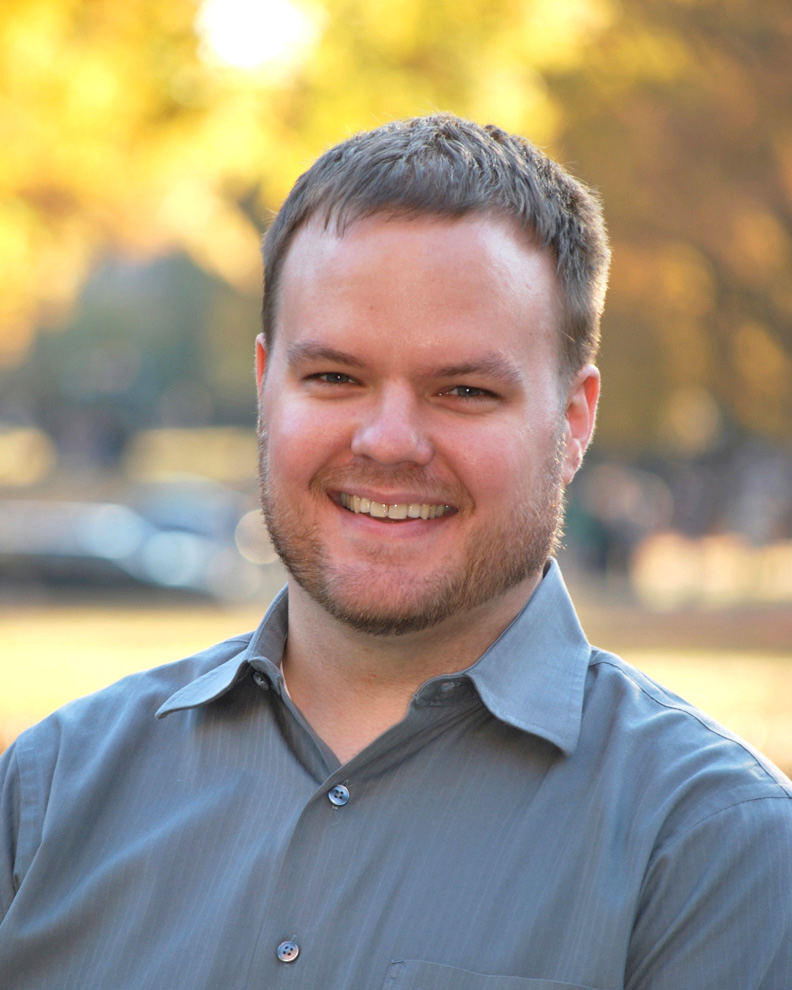}}]{Zachary T. Harmany}
Zachary Harmany received the B.S. (magna cum lade, with honors) in Electrical Engineering and B.S. (cum lade) in Physics from The Pennsylvania State University in 2006. Currently, he is a Ph.D. student in the department of Electrical and Computer Engineering at Duke University. In 2010 he was a visiting researcher at The University of California, Merced. His research interests include nonlinear optimization, functional neuroimaging, and signal and image processing with applications in medical imaging, astronomy, and night vision. 
\end{IEEEbiography}

\begin{IEEEbiography}[{\includegraphics[width=1in,height=1.25in,clip,keepaspectratio]{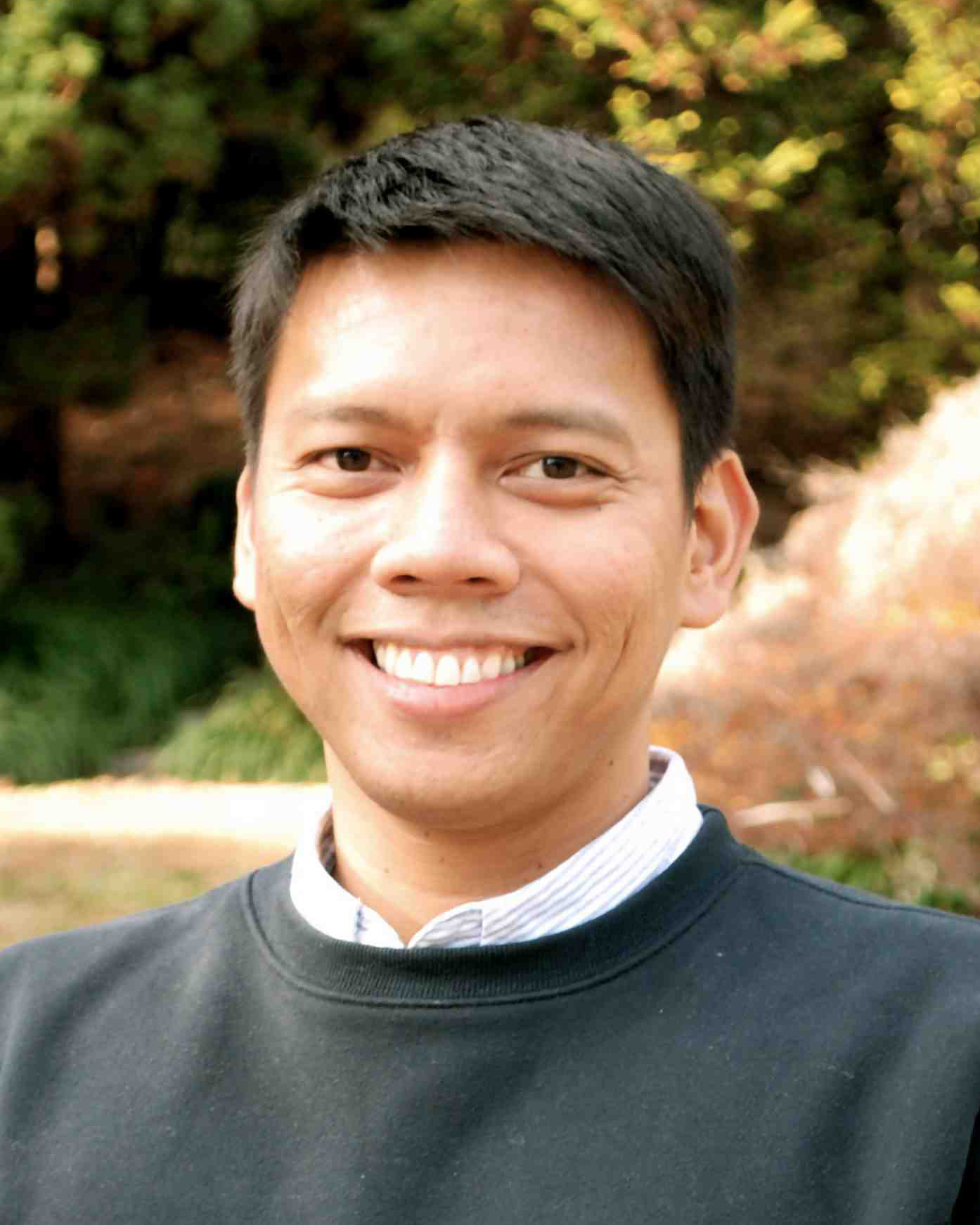}}]{Roummel F. Marcia}
Roummel Marcia received his B.A. in Mathematics from Columbia
University in 1995 and his Ph.D. in Mathematics
from the University of California, San Diego in 2002. He was a
Computation and Informatics in Biology and Medicine
Postdoctoral Fellow in the Biochemistry Department at the University
of Wisconsin-Madison and a Research Scientist
in the Electrical and Computer Engineering at Duke University. He is
currently an Assistant Professor of Applied
Mathematics at the University of California, Merced. His research
interests include nonlinear optimization, numerical
linear algebra, signal and image processing, and computational biology.
\end{IEEEbiography}

\begin{IEEEbiography}[{\includegraphics[width=1in,height=1.25in,clip,keepaspectratio]{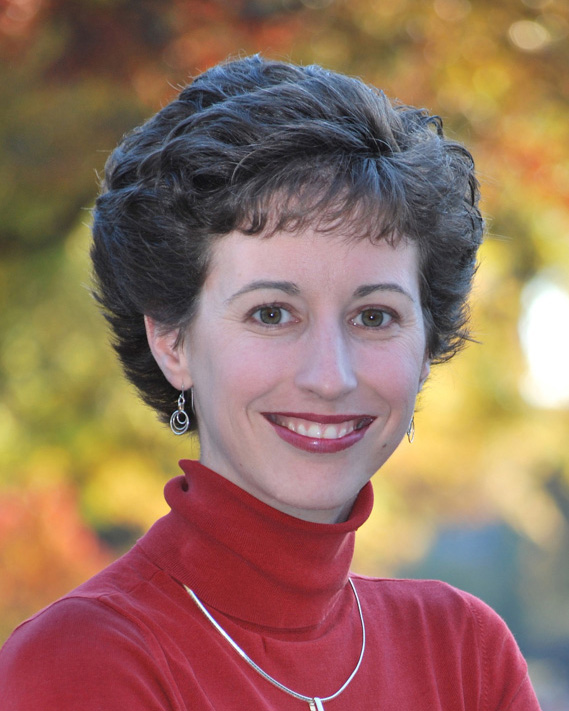}}]{Rebecca M. Willett}
Rebecca Willett is an assistant professor in the Electrical and
Computer Engineering Department at Duke University. She completed her
Ph.D. in Electrical and Computer Engineering at Rice University in 2005.
Prof.~Willett received the National Science Foundation CAREER Award in
2007, is a member of the DARPA Computer Science Study Group, and
received an Air Force Office of Scientific Research Young Investigator
Program award in 2010. Prof.~Willett has also held visiting researcher
positions at the Institute for Pure and Applied Mathematics at UCLA in
2004, the University of Wisconsin-Madison 2003-2005, the French
National Institute for Research in Computer Science and Control
(INRIA) in 2003, and the Applied Science Research and Development
Laboratory at GE Healthcare in 2002. Her research interests include
network and imaging science with applications in medical imaging,
wireless sensor networks, astronomy, and social networks. Additional
information, including publications and software, are available online
at \url{http://www.ee.duke.edu/~willett/}.
\end{IEEEbiography}




\end{document}